\documentclass[journal]{IEEEtran}
\pdfoutput=1

\usepackage{amsfonts}

\usepackage{amssymb}
\usepackage{amsmath}
\usepackage{mathrsfs}
\usepackage{graphics}
\usepackage{epsfig}
\usepackage{epstopdf}
\usepackage{amsthm}
\usepackage{textcomp}
\usepackage{color}
\usepackage{float}
\usepackage{multirow}
\usepackage{diagbox}
\usepackage{bm}
\usepackage{soul} %
\usepackage{slashed}
\usepackage{algorithmic}
\usepackage{algorithm}

\newtheorem{rk}{Remark}
\newtheorem{thm}{Theorem}

\newtheorem{ex}{Experiment}

\ifCLASSINFOpdf
\else
\fi
\usepackage{algorithmic}
\usepackage{algorithm}

\usepackage{array}
\begin{document}
%
\title{Cost-effective company response policy for \\ product co-creation in company-sponsored \\ online community}
%
%
%

\author{Jiamin Hu,
		Lu-Xing Yang,~\IEEEmembership{Member,~IEEE,}
        Xiaofan Yang,~\IEEEmembership{Member,~IEEE,}
        Kaifan Huang,
        \\ Gang Li, ~\IEEEmembership{Senior Member,~IEEE,}
        Yong Xiang, ~\IEEEmembership{Senior Member,~IEEE}

\thanks{This work was supported by National Key Research and Development Project of China ( Grant No. 2019YFB1706101) and Natural Science Foundation of China (Grant No. 61572006).}
\thanks{J. Hu, X. Yang, and K. Huang are with the School of Big Data \& Software Engineering, Chongqing University, Chongqing, 400044, China, e-mails: 1639169307@qq.com, xfyang1964@cqu.edu.cn, dr.kaifanhuang@gmail.com.}
\thanks{L.-X. Yang, G. Li, and Y. Xiang are with the School of Information Technology, Deakin University, Melbourne, VIC 3125, Australia, e-mails: y.luxing@deakin.edu.au, ligang@ieee.org, yong.xiang@deakin.edu.au.}
}

\maketitle


\begin{abstract}
Product co-creation based on company-sponsored online community has come to be a paradigm of developing new products collaboratively with customers. In such a product co-creation campaign, the sponsoring company needs to interact intensively with active community members about the design scheme of the product. We call the collection of the rates of the company's response to active community members at all time in the co-creation campaign as a company response policy (CRP). This paper addresses the problem of finding a cost-effective CRP (the CRP problem). First, we introduce a novel community state evolutionary model and, thereby, establish an optimal control model for the CRP problem (the CRP model). Second, based on the optimality system for the CRP model, we present an iterative algorithm for solving the CRP model (the CRP algorithm). Thirdly, through extensive numerical experiments, we conclude that the CRP algorithm converges and the resulting CRP exhibits excellent cost benefit. Consequently, we recommend the resulting CRP to companies that embrace product co-creation. Next, we discuss how to implement the resulting CRP. Finally, we investigate the effect of some factors on the cost benefit of the resulting CRP. To our knowledge, this work is the first attempt to study value co-creation through optimal control theoretic approach. 
\end{abstract}

\begin{IEEEkeywords}
product co-creation, company-sponsored online community, company response policy, cost benefit, state evolutionary model, optimal control model, optimality system, iterative algorithm, convergence
\end{IEEEkeywords}

%
\IEEEpeerreviewmaketitle

\section{Introduction}
%
\IEEEPARstart{T}{he} unprecedented advancement of computer and communications techniques has radically transformed all aspects of human life. In particular, the customers of a company like to share their ideas and experiences about their purchased products through social networking sites, which has a significant influence on the company's reputation and sales. Under the great pressure to survive and grow in the fiercely competitive market, companies can no longer design and develop products exclusively on their own. Instead, they need to interact with customers to co-create products \cite{Prahalad2004a, Prahalad2004b}. In practice, a well-managed product co-creation process can yield successful innovations and substantial cost reduction \cite{Ramaswamy2010, Ramaswamy2014}.   

Company-sponsored online communities are characterized by fluidity and voluntariness \cite{Faraj2011, Lusch2015}. First, the participants of a company-sponsored online community may come and go at their will. Second, the sponsoring company takes responsibility of coordinating value co-creation but lacks authority to issue command \cite{Nambisan2017}. Company-sponsored online communities have come to be an effective initiative to co-create value \cite{YanJ2018}. In the past decade, many real-world company-sponsored online communities, ranging from Dell IdeaStorm \cite{Gangi2009} and Starbucks \cite{Wong2016} to SAP \cite{Tavakoli2017} and Aston Martin Community \cite{Essamri2019}, have sprung up.

\subsection{Problem formulation}

Suppose a company intends to co-create a product with the participants of its sponsored online community. At the beginning, the company delivers a product design draft to the community. Then, the company interacts intensively with the community participants in the following iterative way: some participants make valuable suggestions about the product design scheme, the company responds by drawing on the suggestions to revise the scheme, and the like. The interaction process continues until the product co-creation activity terminates. As a result, a perfect product design scheme comes into being \cite{Bhalla2011}. 

In the above product co-creation process, the company's continuous response helps to boost the participants' engagement enthusiasm and hence yield a perfect product design scheme. In the context, we call the collection of the company's response rates at all time of the product co-creation activity as a \emph{company response policy} (CRP). We dream of an optimal CRP, i.e., a CRP with the highest cost benefit. However, the dream never comes true. This is because there are so many candidate CRPs that finding an optimal CRP from these candidate CRPs is computationally prohibitive. Therefore, we are forced to take a step back. Specifically, by defining a satisfactory CRP as a CRP that outperforms most other candidate CRPs in terms of the cost benefit, we are satisfied with dealing with the following weaker problem: 

\emph{Company response policy (CRP) problem: Finding a satisfactory CRP from all the candidate CRPs}.

This is a novel problem. This paper aims to address the problem.

\subsection{Contributions}

The main contributions of this work are outlined below.

\begin{itemize}
	\item The CRP problem is reduced to an optimal control model. First, by using the epidemic modeling technique and taking the effect of CRP into account, we introduce a novel `epidemic' model for capturing the evolutionary process of state of the online community. Second, based on the introduced `epidemic' model, we estimate the cost benefit of a CRP. Finally, with the intent of finding a CRP with the highest cost benefit, we establish an optimal control model for the CRP problem (the CRP model).
	\item The CRP model is solved. First, by using the well-known Pontryagin Maximum Principle, we derive the optimality system for the CRP model. Second, by invoking a standard procedure for solving optimal control problems, we present an iterative algorithm for solving the CRP model. Next, through extensive numerical experiments, we corroborate the convergence and effectiveness of the presented algorithm. Therefore, we recommend the resulting CRP to those companies that embrace product co-creation. Finally, we discuss how to implement the resulting CRP. Additionally, we examine the effect of some factors on the cost benefit of the resulting CRP. 
\end{itemize}

The remaining materials are organized in this fashion: Section II reviews the related work. Section III establishes the CRP model. In section IV, the optimality system for the CRP model is derived, and an algorithm for solving the CRP model is presented. Section V validates the convergence and effectiveness of the algorithm. Section VI makes further discussions. Section VII closes this work.

\section{Related work}

This section aims to review the related work and highlight the innovations of the present paper.

\subsection{Online co-creation communities}

Online co-creation communities (O3Cs, for short) can be classified as two categories: \emph{autonomous} and \emph{company-sponsored} \cite{Zwass2010}. The participants of an autonomous O3C (AO3C, for short) perform value co-creation activities independently of any companies. In contrast, value co-creation activities in a company-sponsored O3C (CSO3C, for short) are directed for the benefit of the sponsoring company. 

AO3Cs have received considerable interests from value co-creation researchers and practitioners. To name just a few examples, \cite{Porter2008} suggested a method for cultivating the trust of an AO3C in a company, \cite{Faraj2015} explored the formation mechanism of collaboration in an AO3C, \cite{Constantinides2015, Fernandes2016} proposed some effective measures to enhance the engagement of the participants of an AO3C in value co-creation, and \cite{Spaulding2010, Fisher2019, Rubio2020} identified some key factors that help to co-create value in an AO3C. To our knowledge, the previous researches on AO3Cs are all empirical and the resulting conclusions are all qualitative. This indicates that our understanding of AO3Cs is still in its infancy. 

Also, CSO3Cs have attracted intense attentions. To name but a few examples, \cite{Jeppesen2006, Franke2013} revealed some motives for an individual to participate in a CSO3C, \cite{Manchanda2015} presented an empirical evidence for the effectiveness of value co-creation in a CSO3C, \cite{YanJ2018} validated that company's internal employees play an important role in CSO3C-based value co-creation. \cite{Priharsari2020} made a comprehensive survey on the references that are related to CSO3C and were published before 2020. In particular, \cite{ChenL2014} is closely related to the present paper. In this paper, the authors empirically found that the response rate of the sponsoring company to the participants' valuable suggestions has significant influence on both the community contribution level and the duration of active participation. Still, the extant efforts on CSO3Cs are all empirical, indicating that the study of CSO3Cs is in the early stage. 

Mathematics is the foundation of all sciences. The science of value co-creation is no exception. In our opinion, the very key to revealing the mystry of value co-creation lies in establishing and studying mathematical models characterizing the major features of value co-creation.

\subsection{Epidemic modeling}

Mathematical analysis of epidemic diseases originates from the work by Daniel Bernoulli, the famous Swiss mathematician, about the effectiveness evaluation of a then popular inoculation procedure against the smallpox virus \cite{Dietz2002}. At the beginning of the 20th century, differential and difference equations started to be used to model and analyze the spread of epidemics. In particular, in the seminal work by Kermack and McKendrick \cite{Kermack1927, Kermack1933}, the authors suggested the classical susceptible-infected-removed (SIR) model for the spread of epidemics with no recurrence and the classical susceptible-infected-susceptible (SIS) model for the spread of recurrent and endemic diseases, and established the threshold theory about the SIS model. This work laid a solid foundation for the mathematical theory of epidemic modeling. 

From then on, the classical epidemic models have been extended and generalized rapidly toward different directions, e.g., from population-based epidemic models to network-based epidemic models \cite{Keeling2005, Pare2020}, from deterministic epidemic models to stochastic epidemic models \cite{Allen2008, Britton2010}, from spread of epidemics to resource allocation for the control of epidemics \cite{Nowzari2017, Somers2021}, and from epidemic models on static networks to epidemic models on dynamic networks \cite{Pare2018, Leitch2019}. See \cite{Zino2022} for a comprehensive survey on the theory of epidemic modeling. As one of the main goals of epidemic modeling, different epidemic models have been developed for understanding and containing the spread of various real-life epidemic diseases such as SARS \cite{Hsieh2004} and COVID-19 \cite{Gatto2020}.    

In nature, human society, and engineering area, there exist a wide spectrum of propagation phenomena other than the spread of biological epidemics, ranging from the diffusion of information \cite{JiangY2015} to the propagation of failure \cite{XingL2013}. In the elegant paper published in 1964 by Goffman and Newill \cite{Goffman1964}, an epidemic process is defined as transition from one state to another owing to exposure to some phenomenon. For example, the authors compared the transmission of ideas to the spread of infectious diseases. Finally, the authors elucidated the importance of inspecting various propagation processes through epidemic modeling. Since then, a multitude of `epidemic' models, ranging from rumor spreading models \cite{Daley1964, YangLX2018a} and malware propagation models \cite{Kephart1991, YuS2015} to word-of-mouth propagation models \cite{LiP2018a, ZhangTR2019}, have emerged.  

Every participant of a CSO3C must be in one of two possible states: \emph{active}, i.e., actively engaging in the design of the product, and \emph{inactive}, i.e., not active. We define the state of the community as the collection of the states of all participants of the community. Following \cite{Goffman1964}, the process of transition between the active state and the inactive state can be viewed as an `epidemic' process and hence can be studied through epidemic modeling. To our knowledge, to date there is no such an `epidemic' model.

\subsection{Optimal control theory combined with epidemic modeling}

Optimal control theory as an integral part of optimization theory is dedicated to finding a time-varying control strategy of a dynamical system so that a given performance index is optimized \cite{Liberzon2012}. Optimal control theory when combined with epidemic modeling is especially suited to the development of cost-effective control policies of various `epidemic' processes, such as the spread of infectious diseases \cite{Aldila2020, Sameni2022}, the propagation of malware \cite{Eshghi2016, YangLX2021a}, and the diffusion of rumors \cite{ZhaoJ2019, LinY2021}. 

Differential game theory as an extension of optimal control theory studies time-varying strategic interactions between informed and rational players \cite{Basar2012}. Differential game theory when combined with epidemic modeling provides a powerful tool for the design of cost-effective control policies of `epidemic' processes in the presence of strategic adversary, ranging from clarification of rumors in the presence of strategic rumormonger \cite{HuangDW2021} to defense against cyberattacks in the presence of strategic cyber attacker \cite{YangLX2019a, YangLX2021b}.

The prerequisite for developing a cost-effective control policy of an `epidemic' process through game theoretic approach lies in having the ability to acquire quite a number of adversary-related parameters. In practice, these parameters are often unavailable. In contrast, the optimal control theoretic approach to the design of cost-effective control policy of an `epidemic' process is feasible if a few relevant parameters are known, and these parameters may be estimated relatively accurately by leveraging relevant historical data. Consequently, we choose to deal with the CRP problem through optimal control approach.

\subsection{Innovations of the present paper}

In the context of product co-creation in company-sponsored online community, a new problem (i.e., the CRP problem) is proposed and studied, with the goal of enhancing the cost benefit of product co-creation. This is the core innovation of this paper. The key step to the solution of the CRP problem is to establish an optimal control model for the CRP problem. For this purpose, the benefit of a CRP needs to be estimated. To this end and in view of the mobility and voluntariness of the product co-creation participants, a new `epidemic' model for characterizing the evolutionary process of state of the company-sponsored online community is introduced. This is another innovation of the present paper. Although there exist similar known epidemic models, our epidemic model has special meaning in the context of product co-creation.

\section{Modeling the CRP problem}

This section is devoted to modeling the CRP problem proposed in section I. First, a CRP is formalized. Second, the state evolutionary model for the company-sponsored online community is proposed. Finally, the CRP problem is modeled as an optimal control problem. See Table I for a list of major notations used in this article and their meanings.

\begin{table}
	\centering
	\caption{Notations and their meanings}
	\begin{tabular}{|c|c|}
		\hline
		$T$ & product co-creation period \\
		\hline
		$x$ & company response policy (CRP) \\
		\hline
		$\overline{x}$ & maximal response rate \\
		\hline
		$X_{T, \overline{x}}$ & set of feasible CRPs \\
		\hline
		$\omega_1$ & standard cost \\
		\hline
		$C(x)$ & overall cost for implementing the CRP $x$ \\
		\hline
		$A(t)$ & number of active participants at time $t$ \\
		\hline
		$I(t)$ & number of inactive participants at time $t$ \\
		\hline
		$\mu$ & community inflow rate \\
		\hline
		$\delta_1$ & first community outflow rate \\
		\hline
		$\delta_2$ & second community outflow rate \\
		\hline
		$\beta_1$ & first influence function \\
		\hline
		$\beta_2$ & second influence function \\
		\hline
		$\alpha$ & inaction rate \\
		\hline
		$\omega_2$ & standard benefit \\
		\hline
		$J(x)$ & cost benefit of implementing the CRP $x$ \\
		\hline
		$H(A, I, x, \lambda_1, \lambda_2)$ & Hamiltonian function for the CRP model \\
		\hline
		$(\lambda_1, \lambda_2)$ & ajoint for $H$ \\
		\hline
	\end{tabular} 
\end{table}

\subsection{Company response policy}

Consider the CRP problem. Suppose the company intends to organize a product co-creation activity in its sponsored online community. At the beginning, the company delivers a product design draft to the community. Next, the company interacts with the community participants in the following iterative way: some participants propose suggestions for the revision of the product design scheme, the company responds by drawing on good suggestions to revise the scheme, and the like. The interaction process continues until the activity terminates. 

Suppose the above-mentioned product co-creation activity starts from the initial time $t = 0$ and terminates at the time $t = T$. We call the parameter $T$ as the \emph{product co-creation period}. For $0 \leq t \leq T$, let $x(t)$ denote the company response rate at time $t$. We call the function $x$ as a \emph{company response policy} (\emph{CRP}, for short). 

On the one hand, for ease in implementation, assume a feasible CRP is piecewise continuous. Let $PC[0, T]$ denote the set of piecewise continuous functions defined on $[0, T]$. On the other hand, owing to the limited product co-creation budget, assume a feasible CRP is bounded from above. Let $\overline{x}$ denote the common upper bound on the response rates of a feasible CRP at all time of the product co-creation process. We call the parameter $\overline{x}$ as the \emph{maximal response rate}. In a word, a CRP is feasible if and only if it is in the set
\begin{equation}
	X_{T, \overline{x}} = \{x \in PC[0, T] : 0 \leq x(t) \leq \overline{x}, 0 \leq t \leq T\}.
\end{equation}

The implementation of a feasible CRP comes at a cost. Assume the rate of increase of the cost at any time is proportional to the company's response rate at that time. Then the total cost for implementing the CRP $x$ equals 
\begin{equation}
	C(x) = \omega_1\int_0^T x(t)dt. 
\end{equation}
Here, $\omega_1$ represents the per-unit-time cost for the company response, which is assumed to be positive and constant. We call the parameter $\omega_1$ as the \emph{standard cost}.

\subsection{A state evolutionary model for the online community}

At any time of the product co-creation activity, each participant of the online community is in one of two possible states: \emph{active}, i.e., actively participating in the revision of the product design scheme by proposing valuable suggestions in a time interval of a given small length, and \emph{inactive}, i.e., not active in a time interval of the same length. Let $A(t)$ and $I(t)$ denote the number of active and inactive participants at time $t$, respectively. We call the ordered pair $(A(t), I(t))$ as the \emph{state of the community at time $t$}. The initial community state, $(A(0), I(0))$, can be determined by a company worker through continuously observing the active statuses of all community participants in a time interval preceding the initial time $t= 0$.

On the one hand, since any individual may enter or exit the community  at his or her will, $A(t)$ and $I(t)$ vary over time. On the other hand, the state of every participant of the community is varying over time. Specifically, an active participant may become inactive owing to various reasons (heavy work, being on vacation, being tired, etc.), and an inactive participant may become active owing to the influence of either the company response or the active participants. As a result, $A(t)$ and $I(t)$ are varying over time as well. 

An epidemic model is established following a two-step procedure: First, introduce a set of epidemic-related assumptions and parameters. Second, use these assumptions and parameters to describe an epidemic model \cite{Zino2022}. For the purpose of establishing an epidemic model for characterizing the evolutionary process of the state of the community over time, let us introduce a set of assumptions and parameters as follows. 

\begin{enumerate}
\item[(i)] Assume an outside individual has no knowledge of the product design scheme when entering the community. As a result, the individual is inactive at that time. 
\item[(ii)] Let $\mu$ denote the per-unit-time number of outside individuals entering the community, which is assumed to be positive and constant. We refer to $\mu$ as the \emph{community inflow rate}.
\item[(iii)] Let $\delta_1$ denote the per-unit-time probability of an active community participant exiting the community, which is assumed to be positive and constant. We refer to $\delta_1$ as the \emph{first community outflow rate}.
\item[(iv)] Let $\delta_2$ denote the per-unit-time probability of an inactive community participant exiting the community, which is assumed to be positive and constant. We refer to $\delta_2$ as the \emph{second community outflow rate}. Intuitively, $\delta_2 > \delta_1$.
\item[(v)] Let $\beta_1(z)$ denote the per-unit-time probability of an inactive participant becoming active owing to the influence of the company response rate $z$. We refer to the function $\beta_1$ as the \emph{first influence function}. Intuitively, $\beta_1$ is monotonically increasing and flattens out.
\item[(vi)] Let $\beta_2(z)$ denote the per-unit-time probability of an inactive participant becoming active owing to the influence of the number $z$ of active participants. We call the function $\beta_2$ as the \emph{second influence function}. Intuitively, $\beta_2$ is monotonically increasing and flattens out. For technical reasons, assume $\beta_2$ is differentiable. 
\item[(vii)] Let $\alpha$ denote the per-unit time probability of an active participant becoming inactive owing to various reasons, which is assumed to be positive and constant. We refer to $\delta$ as the \emph{inaction rate}.
\end{enumerate}

On the one hand, the assumptions (ii)-(iv) quantitatively characterize the fluidity of the company-sponsored online community. On the other hand, the assumptions (v)-(vii) quantitatively characterize the voluntariness of a community participant to be active or inactive. See references \cite{Faraj2011, Lusch2015, Nambisan2017}. This set of assumptions provides a solid foundation for the subsequent modeling of the CRP problem.

\begin{rk}
  In practice, the inflow rate, the two outflow rates, and the inaction rate are all time-varying. In our modeling, the four rates are averaged over time and hence are assumed to be constant. Anyway, averaging is a commonly used technique in establishing epidemic models.
\end{rk}

The above assumptions and parameters imply the following result.

\begin{thm}
Under the influence of the CRP $x$, the community state evolves obeying the following rule:
\begin{equation}
\left\{\begin{aligned}
\frac{dA(t)}{dt} &= \left[\beta_1(x(t)) + \beta_2(A(t)) \right]I(t) - \alpha A(t) - \delta_1 A(t),\\
\frac{dI(t)}{dt} &= \mu - \left[\beta_1(x(t)) + \beta_2(A(t)) \right]I(t) + \alpha A(t) - \delta_2 I(t), \\
& 0 \leq t \leq T, \\
A(0) &= A_0, I(0) = I_0.
\end{aligned}\right.
\end{equation}
\end{thm}

We call the system (3) as the \emph{community state evolutionary model}. See Fig. 1 for a diagram of the model. According to Goffman and Newill's opinion \cite{Goffman1964}, this is an `epidemic' model, because it characterizes the `epidemic' process of transitions between the active state and the inactive state. According to the taxonomy of epidemic models proposed in \cite{Zino2022}, the state evolutionary model of the online community is a population-based `epidemic' model. The main function of the community state evolutionary model (3) is to accurately predict the community state in the future time interval $[0, T]$. Without the model, it would be difficult to accomplish the goal. As any long-term prediction lacks prediction accuracy, the parameter $T$ should take on a relatively small value to guarantee a relatively high prediction accuracy of the community state evolutionary model. 

\begin{figure}[!ht]
	\setlength{\abovecaptionskip}{0.cm}
	\setlength{\belowcaptionskip}{-0.cm}
	\centering
	\includegraphics[width=0.38\textwidth]{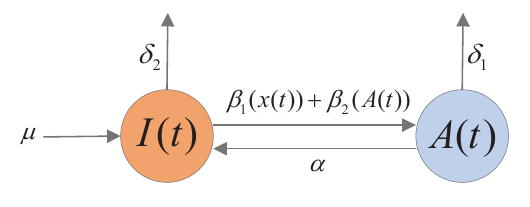}
	\caption{Diagram of the community state evolutionary model.}
\end{figure}

\subsection{Optimal control model for the CRP problem}

In the CRP problem, all active participants make contributions to the final product design scheme. Let $\omega_2$ be the per-unit-time benefit brought by an active participant, which is assumed to be positive and constant. We call the parameter $\omega_2$ as the \emph{standard benefit}. The assumption that the per-unit-time contribution is identical for all active participants implies the per-unit-time contribution is averaged over all active participants. Although this assumption is oversimplified, it is trivial to generalize the yielded model to the more general situation where different active participants may make different per-unit-time contributions. Based on the accurate prediction of the community state evolutionary model (3), the cost benefit of implementing the CRP $x$ is expected to be
\begin{equation}
  J(x) = \omega_2\int_0^T A(t)dt - \omega_1\int_0^T x(t)dt.
\end{equation}
Therefore, the CRP problem may be reduced to the following open-loop, deterministic optimal control problem:
\begin{equation}
\begin{aligned}
& \max_{x \in X_{T, \overline{x}}} J(x) = \omega_2\int_0^T A(t)dt - \omega_1 \int_0^T x(t)dt \\
& s. t. \left\{
\begin{aligned}
\frac{dA(t)}{dt} &= \left[\beta_1(x(t)) + \beta_2(A(t)) \right]I(t) - \alpha A(t) - \delta_1 A(t),\\
\frac{dI(t)}{dt} &= \mu - \left[\beta_1(x(t)) + \beta_2(A(t)) \right]I(t) + \alpha A(t) \\
& \quad - \delta_2 I(t), \quad 0 \leq t \leq T, \\
A(0) &= A_0, I(0) = I_0.
\end{aligned}
\right.
\end{aligned}
\end{equation}

We call the optimal control problem (5) as the \emph{CRP model}. Every instance of the CRP model (CRP instance, for short) is represented by a 12-tuple of the form
\begin{equation}
	\mathbb{M} = (A_0, I_0, T, \overline{x}, \mu, \delta_1, \delta_2, \alpha, \beta_1, \beta_2, \omega_1, \omega_2).
\end{equation}

For this purpose, the prediction accuracy of the community state evolutionary model (3) must be guaranteed. At the first sight, the aim of the CRP model (5) is to find an optimal control function, i.e., a control function that maximizes the objective functional. But this is not the case. First, since the CRP model is inherently complex and there are numerous feasible control functions, finding an optimal control function is computationally prohibitive. Second, since the aim of the CRP problem is to find a satisfactory CRP, and the CRP model is a mathematical model for the CRP problem, we should be satisfied with finding a satisfactory control function, which stands for a satisfactory CRP. For this purpose, the community state evolutionary model (3) must own a guaranteed prediction accuracy.

\section{Solving the CRP model}

In the preceding section, we established the CRP model. This section is committed to solving the CRP model. 

The standard procedure of solving an open-loop, deterministic optimal control problem is as follows: First, derive the optimality system for the optimal control problem. Second, solve the optimality system. See \cite{Liberzon2012}. The CRP model (5) is an open-loop, deterministic optimal control problem. Below let us follow the procedure to solve the CRP model.

\subsection{The optimality system for the CRP model}

The Hamiltonian function for the CRP model (5) reads
\begin{equation}
\begin{aligned}
& \quad H(A, I, x, \lambda_1, \lambda_2) \\
& = \omega_2 A - \omega_1 x + \lambda_1\left\{\left[\beta_1(x) + \beta_2(A)\right]I - \alpha A - \delta_1 A\right\} \\
& \quad + \lambda_2 \left\{\mu - \left[\beta_1(x) + \beta_2(A)\right]I + \alpha A - \delta_2 I\right\}.
\end{aligned}
\end{equation}
Here, $(\lambda_1, \lambda_2)$ is the adjoint variable for $H$.

We have the following result.

\begin{thm}
Let $x$ be an optimal CRP for the CRP model (5). Let $(A, I)$ be the solution to the associated community state evolutionary model (3). Then, there exists an adjoint function $(\lambda_1, \lambda_2)$ such that
\begin{equation}
\left\{\begin{aligned}
& \frac{d\lambda_1(t)}{dt} = -\omega_2 + [\alpha + \delta_1 - \beta_2^{'}(A(t))I(t)]\lambda_1(t) \\
& \quad - [\alpha - \beta_2^{'}(A(t))I(t)]\lambda_2(t), \quad 0 \leq t \leq T, \\
&\frac{d\lambda_2(t)}{dt} = -[\beta_1(x(t)) + \beta_2(A(t))]\lambda_1(t) \\
& \quad + [\delta_2 + \beta_1(x(t)) + \beta_2(A(t))]\lambda_2(t), \quad 0 \leq t \leq T, \\
&\lambda_1(T) = \lambda_2(T) = 0.
\end{aligned}\right.
\end{equation}
Moreover,
\begin{equation}
\left\{\begin{aligned}
x(t) & \in\arg\max\limits_{0 \leq \tilde{x} \leq \overline{x}} \left\{[\lambda_1(t) - \lambda_2(t)]I(t)\beta_1(\tilde{x}) - \omega_1 \tilde{x}\right\}, \\
& 0 \leq t \leq T.
\end{aligned}\right.
\end{equation}
\end{thm}

\begin{IEEEproof}
First, $\lambda_1(T) = \lambda_2(T) = 0$ follows from the unspecified terminal cost and the free final state. Second, the Pontryagin Maximum Principle \cite{Liberzon2012} implies there is an adjoint function $(\lambda_1, \lambda_2)$ such that
\begin{equation}
\left\{\begin{aligned}
\frac{d\lambda_1(t)}{dt} &= -H_A(A(t), I(t), x(t), \lambda_1(t), \lambda_2(t)), \\
\frac{d\lambda_2(t)}{dt} &= -H_I(A(t), I(t), x(t), \lambda_1(t), \lambda_2(t)), \\
& 0 \leq t \leq T.
\end{aligned}\right.
\end{equation}
The first two equations in the system (8) follow by direct calculations. Finally, the Pontryagin Maximum Principle implies
\begin{equation}
	x \in \arg\max\limits_{\tilde{x} \in X}H(A, I, \tilde{x}, \lambda_1, \lambda_2),
\end{equation}
The system (9) follows through simple calculations.
\end{IEEEproof}

According to Theorem 2, the optimality system for the CRP model (5) is shown in Eqs. (12). The optimality system may be viewed as a system in the CRP $x$, where both the state function $(A, I)$ and the adjoint function $(\lambda_1, \lambda_2)$ play the role of auxiliary function.
\begin{figure*}
\begin{equation}
  \left\{
  \begin{aligned}
  \frac{dA(t)}{dt} &= \left[\beta_1(x(t)) + \beta_2(A(t)) \right]I(t) - \alpha A(t) - \delta_1 A(t), \quad 0 \leq t \leq T, \\
  \frac{dI(t)}{dt} &= \mu - \left[\beta_1(x(t)) + \beta_2(A(t)) \right]I(t) + \alpha A(t) - \delta_2 I(t), \quad 0 \leq t \leq T, \\
\frac{d\lambda_1(t)}{dt} &= -\omega_2 + [\alpha + \delta_1 - \beta_2^{'}(A(t))I(t)]\lambda_1(t) - [\alpha - \beta_2^{'}(A(t))I(t)]\lambda_2(t), \quad 0 \leq t \leq T, \\
\frac{d\lambda_2(t)}{dt} &= -[\beta_1(x(t)) + \beta_2(A(t))]\lambda_1(t) + [\delta + \beta_1(x(t)) + \beta_2(A(t))]\lambda_2(t), \quad 0 \leq t \leq T, \\
x(t) &\in\arg\max\limits_{0 \leq \tilde{x} \leq \overline{x}} \left\{[\lambda_1(t) - \lambda_2(t)]I(t)\beta_1(\tilde{x}) - \omega_1 \tilde{x}\right\}, \quad 0 \leq t \leq T, \\
A(0) &= A_0, I(0) = I_0, \lambda_1(T) = \lambda_2(T) = 0.
  \end{aligned}
  \right.
\end{equation}
\hrulefill
\end{figure*}

\subsection{An algorithm for solving the CRP model}

It follows from Theorem 2 that an optimal CRP for the CRP model (5) is a solution to the optimality system (12). In the situation where the optimality system admits a sole solution, this solution is an optimal CRP. However, if the optimality system admits more than one solution, a solution to the optimality system is not necessarily an optimal CRP for the CRP model. Owing to the inherent complexity of the optimality system, any attempt to prove the uniqueness of solution to this system proves futile. As a matter of fact, the system may admit multiple solutions. As a result, the CRP model is not equivalent to the optimality system. Nevertheless, the task of solving the CRP model may be reduced to two successive subtasks: solving the optimality system and checking the cost effectiveness of the yielded solution.

Owing to the complexity of the optimality system (12), any attempt to analytically solve the system proves futile. Hence, we turn to numerically solve the system. Invoking the Forward-Backward Sweep Method for solving the optimality system for an optimal control problem \cite{McAsey2012}, we design an algorithm for numerically solving the optimality system (12). The basic idea of our algorithm is to yield a sequence of CRPs. In the case where the sequence converges, the finally yielded CRP is a numerical solution to the optimality system (12). In algorithm 1 the algorithm is sketched and is referred to as the CRP algorithm.  

\begin{algorithm}[h] 
	\caption{CRP}   
	\label{alg1}   
	{\bf Input:} an instance $\mathbb{M} = (A_0, I_0, T, \overline{x}, \mu, \delta_1, \delta_2, \alpha, \beta_1, \beta_2, \omega_1,$ $\omega_2)$ of the CRP model (5), convergence error $\epsilon$. \\
	{\bf Output:} a control policy $x$.	
	\begin{algorithmic}[1] 
		\STATE $k \leftarrow 0$; $x^{(0)} \leftarrow 0$;
		\REPEAT 
		\STATE $k \leftarrow k + 1$;
		\STATE  forwardly calculate the state function $(A, I)$ using Eqs. (3) with $x \leftarrow x^{(k-1)}$; $(A^{(k)}, I^{(k)}) \leftarrow (A, I)$;
		\STATE backwardly calculate the adjoint function $(\lambda_1, \lambda_2)$ using Eqs. (8) with $x \leftarrow x^{(k-1)}$ and $(A, I) \leftarrow (A^{(k)}, I^{(k)})$; $(\lambda_1^{(k)}, \lambda_2^{(k)}) \leftarrow (\lambda_1, \lambda_2)$;
		\STATE calculate the control policy $x$ using Eqs. (9) with $(A, I) \leftarrow (A^{(k)}, I^{(k)})$, and $(\lambda_1, \lambda_2) \leftarrow (\lambda_1^{(k)}, \lambda_2^{(k)})$; $x^{(k)} \leftarrow x$;
		\UNTIL{$\sup_{0 \leq t \leq T}|x^{(k)}(t) - x^{(k - 1)}(t)| < \epsilon$}; 		
		\RETURN $x^{(k)}$.
	\end{algorithmic} 
\end{algorithm} 

The overall time cost of the CRP algorithm is dominated by that for solving the collection of maximization problems involved in the optimality system (12). Hence, the key to reducing the time cost of the CRP algorithm is to reduce the time cost for solving each of these maximization problems. To this end, we present the following theorem. See \cite{Sethi2000, Burden2015}.
  
\begin{thm}
Let $x$ be a solution to the optimality system (12), $(A, I)$ the associated state function, $(\lambda_1, \lambda_2)$ the associated adjoint function. For each $t \in [0, T]$, the following claims hold.
\begin{itemize}
  \item[(i)] If $[\lambda_1(t) - \lambda_2(t)]I(t)\beta_1^{'}(\overline{x}) > \omega_1$, then $x(t) = \overline{x}$.
  \item[(ii)] If $[\lambda_1(t) - \lambda_2(t)]I(t)\beta_1^{'}(0) < \omega_1$, then $x(t) = 0$.
  \item[(iii)] If $[\lambda_1(t) - \lambda_2(t)]I(t)\beta_1^{'}(\overline{x}) < \omega_1$,
  $[\lambda_1(t) - \lambda_2(t)]I(t)\beta_1^{'}(0)$ $> \omega_1$, then $x(t) = \left(\beta_1^{'}\right)^{-1}\left(\frac{\omega_1}{(\lambda_1(t) - \lambda_2(t))I(t)}\right)$.
\end{itemize}
\end{thm}

\begin{IEEEproof}
  Let $G_t(\tilde{x}) = [\lambda_1(t) - \lambda_2(t)]I(t)\beta_1(\tilde{x}) - \omega_1 \tilde{x}$.

(i) In this case, $G_t^{'}(\tilde{x}) > 0$ for $0 \leq \tilde{x} \leq \overline{x}$. As a result, $G_t(\tilde{x})$ is strictly increasing in the interval $[0, \overline{x}]$. The claimed result follows.

(ii) In this case, $G_t^{'}(\tilde{x}) < 0$ for $0 \leq \tilde{x} \leq \overline{x}$. As a result, $G_t(\tilde{x})$ is strictly decreasing in the interval $[0, \overline{x}]$. The claimed result follows.

(iii) In this case, $G_t(\tilde{x}) > 0$ for $0 \leq \tilde{x} \leq \frac{\omega_1}{(\lambda_1(t) - \lambda_2(t))I(t)}$, and  $G_t(\tilde{x}) < 0$ for $\frac{\omega_1}{(\lambda_1(t) - \lambda_2(t))I(t)} \leq \tilde{x} \leq \overline{x}$. As a result, $G_t(\tilde{x})$ is first strictly increasing then strictly decreasing, and the turning point is $\frac{\omega_1}{(\lambda_1(t) - \lambda_2(t))I(t)}$. The claimed result follows.
\end{IEEEproof}

This theorem helps to greatly reduce the time cost used for solving the maximization problem (9) and hence the time cost for solving the optimality system.

\begin{rk}
There are two fundamentally different approaches to dealing with an actual optimal control problem. The first approach is to reduce the problem to a continuous-type optimal control model, derive a continuous-type optimality system for the model, apply the classical forward-backward sweep method to the optimality system to present a `continuous-type algorithm' for finding a satisfactory control for the model, and finally discretize the `algorithm' to yield a numerical algorithm. The second approach is to reduce the problem to a discrete-type optimal control model, derive a discrete-type optimality system for the model, and apply the forward-backward sweep method to the optimality system to present a numerical algorithm for finding a satisfactory control for the model. Owing to three reasons, in this paper we choose the first approach to deal with the CRP problem. First, this is a commonly used approach to deal with open-loop, deterministic optimal control models. Second, the results of the relevant comparative experiments reported in the subsequent section show that the CRP obtained through the approach is satisfactory in terms of the cost benefit. Finally, a few profound results about the structure of the resulting CRP (i.e., Theorem 3) are derived analytically. By the way, the CRP algorithm comes from a direct application of the forward-backward sweep method.
\end{rk}

\begin{rk}
	When encoding the CRP algorithm, the time interval $[0, T]$ should be split into a collection of sub-intervals, say, $[0, \frac{T}{N}]$, $[\frac{T}{N}, \frac{2T}{N}]$, $\cdots$, $[\frac{(N - 1)T}{N}, T]$, and the calculations involved in the CRP algorithm should be performed at these split points.
\end{rk}

\section{Convergence, effectiveness, and implementation of the CRP algorithm}
In the preceding section, we presented an algorithm for solving the CRP model (i.e., the CRP algorithm). For the algorithm to work properly, the sequence of CRPs yielded by the algorithm must converge, and the finally yielded CRP must be superior to most other feasible CRPs in terms of the cost benefit. Owing to the highly complex structure of the optimality system (12), it is extremely difficult, if not impossible, to rigorously prove the convergence and effectiveness of the CRP algorithm. Instead, this section turns to experimentally inspect the convergence and effectiveness of the algorithm, followed by a discussion of some issues about the implementation of this algorithm.
\subsection{Convergence and effectiveness of the CRP algorithm}

\begin{ex}
Consider the following CRP instance.
\begin{equation*}
\begin{aligned}
  \mathbb{M}_1 &= (50, 10000, 50, 10, 12, 0.0001, 0.001, 0.1, \\
  &0.05\arctan (0.3z), 0.01\ln (0.01z+1), 10^3, 20).
\end{aligned}
\end{equation*}
First, we observe that, when run on $(\mathbb{M}_1, 10^{-6})$, the CRP algorithm converges in four iterative steps. The resulting sequence of CRPs, denoted $\{x^{(k)}\}_{k=1}^4$, is plotted in Fig. 2(a). Let $x_1^* = x^{(4)}$ denote the finally resulting CRP. Second, randomly and uniformly generate 100 feasible CRPs, denoted $x_1, \cdots, x_{100}$. We observe from Fig. 2(b) that $J(x_1^*) > J(x)$, $x \in \{x_1, \cdots, x_{100}\}$. 
\end{ex}

\begin{figure}[!ht]
	\setlength{\abovecaptionskip}{0.cm}
	\setlength{\belowcaptionskip}{-0.cm}
	\centering
	\includegraphics[scale=0.66]{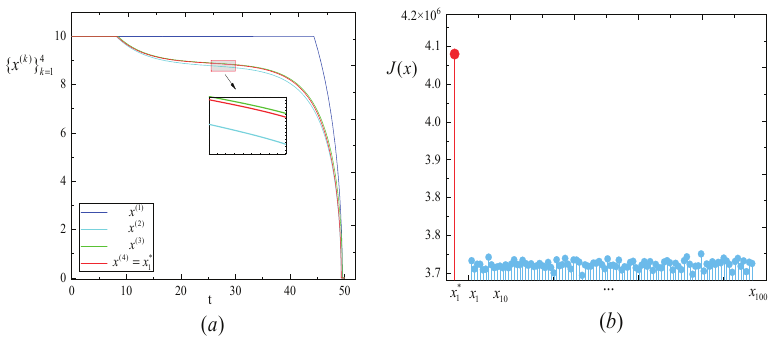}
	\caption{The experimental results in Experiment 1: (a) the resulting sequence of CRPs, denoted $\{x^{(k)}\}_{k=1}^4$, (b) $J(x)$ versus $x$, $x \in \{x_1^*, x_1, \cdots, x_{100}\}$.}
\end{figure}

\begin{ex}
Consider the following CRP instance.
\begin{equation*}
\begin{aligned}
	\mathbb{M}_2 &= (100, 10000, 80, 15, 15, 0.0001, 0.001, 0.15, 0.06z^{1/4}, \\
	&0.003z^{1/3}, 1200, 20).
\end{aligned}
\end{equation*}
First, we observe that, when run on $(\mathbb{M}_2, 10^{-6})$, the CRP algorithm converges in five iterative steps. The resulting sequence of CRPs, denoted $\{x^{(k)}\}_{k=1}^{5}$, is plotted in Fig. 3(a). Let $x_2^* = x^{(5)}$ denote the finally resulting CRP. Second, randomly and uniformly generate 100 feasible CRPs, denoted $x_1, \cdots, x_{100}$. We observe from Fig. 3(b) that $J(x_2^*) > J(x)$, $x \in \{x_1, \cdots, x_{100}\}$.  
\end{ex}

\begin{figure}[!ht]
	\setlength{\abovecaptionskip}{0.cm}
	\setlength{\belowcaptionskip}{-0.cm}
	\centering
	\includegraphics[scale=0.65]{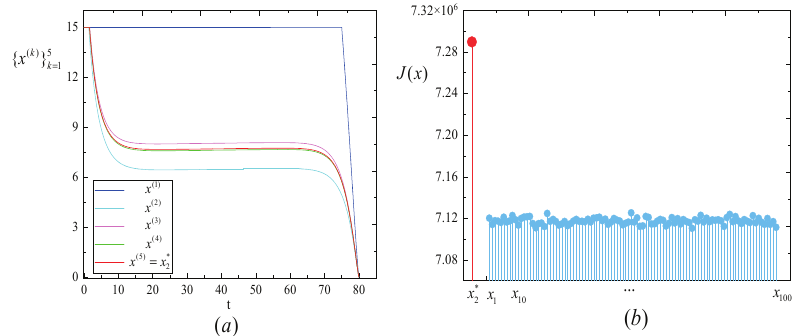}
	\caption{The experimental results in Experiment 2: (a) the resulting sequence of CRPs, denoted $\{x^{(k)}\}_{k=1}^{5}$, (b) $J(x)$ versus $x$, $x \in \{x_2^*, x_1, \cdots, x_{100}\}$.}
\end{figure}

\begin{ex}
Consider the following CRP instance.
\begin{equation*}
\begin{aligned}
  \mathbb{M}_3 &= (150, 10000, 100, 20, 10, 0.0003, 0.001, 0.2, \\
  &0.04\ln(z+1), 0.04\arctan (0.001z), 1000, 25).
\end{aligned}
\end{equation*}
First, we observe that, when run on $(\mathbb{M}_3, 10^{-6})$, the CRP algorithm converges in five iterative steps. The resulting sequence of CRPs, denoted $\{x^{(k)}\}_{k=1}^{5}$, is plotted in Fig. 4(a). Let $x_3^* = x^{(5)}$ denote the finally resulting CRP. Second, randomly and uniformly generate 100 feasible CRPs, denote $x_1, \cdots, x_{100}$. We observe from Fig. 4(b) that $J(x_3^*) > J(x)$, $x \in \{x_1, \cdots, x_{100}\}$.  
\end{ex}

\begin{figure}[!ht]
	\setlength{\abovecaptionskip}{0.cm}
	\setlength{\belowcaptionskip}{-0.cm}
	\centering
	\includegraphics[scale=0.65]{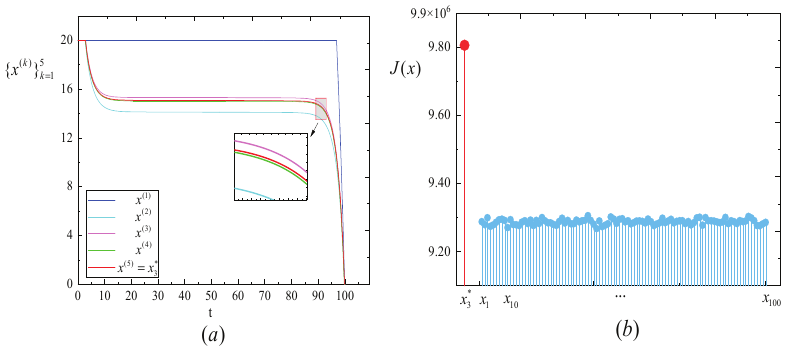}
	\caption{The experimental results in Experiment 3: (a) the resulting sequence of CRPs, denoted $\{x^{(k)}\}_{k=1}^{5}$, (b) $J(x)$ versus $x$, $x \in \{x_3^*, x_1, \cdots, x_{100}\}$.}
\end{figure}

We conducted totally $10000$ similar experiments. In each of these experiments, we observed that (i) the CRP algorithm converges rapidly, and (ii) the finally resulting CRP is superior to 100 randomly generated feasible CRPs. Hence, we conclude that, generally, the CRP algorithm converges rapidly and the resulting CRP is cost-effective. Therefore, we recommend the resulting CRP to companies that embrace product co-creation.

In all of our $10000$ experiments, we observe the phenomenon that the CRP yielded by the CRP algorithm starts at $\overline{x}$ and then decreases. This phenomenon is caused by the fact that the cost benefit yielded by achieving higher response rate at early stage of a product co-creation campaign is higher than that yielded by achieving higher response rate at later stage.

\subsection{Another approach to solving the CRP model}

Dynamic programming provides another method for solving optimal control problems \cite{Bertsekas2012}. The dynamic programming solution of a continuous-type optimal control problem consists of two steps. In the first step, approximate the original problem by a discrete-type optimal control problem. In the second step, recursively solve the discrete-type Hamilton-Jacobi-Bellman equation for this problem to get a discrete optimal control.

For the purpose of using dynamic programming technique to solve the CRP model (5), we need to introduce a set of notations as follows. 

\begin{itemize}
	\item Let $t_0 = 0, t_1 = \frac{T}{N}, t_2 = \frac{2T}{N}, \cdots, t_N = T$ denote the sequence of time points. 
	\item For $0 \leq i \leq N$, let $A_i = A(t_i)$, $I_i = I(t_i)$.
	\item Let $s_0 = 0, s_1 = \frac{S}{M}, s_2 = \frac{2S}{M}, \cdots, s_M = S$ denote the admissible values of all $A_i$ and $I_i$, where $S$ stands for the maximum possible size of the online community.
	\item Let  $p_0 = 0, p_1 = \frac{\overline{x}}{P}, p_2 = \frac{2\overline{x}}{P}, \cdots, p_P = \overline{x}$ denote the admissible values of a control policy $x$ at any time.
	\item For $0 \leq i \leq N$, $0 \leq j, k \leq M$, let $x(i, j, k)$ denote the value of an optimal control policy $x$ at the time point $t_i$ and in the state $(s_j, s_k)$, let $J(i, j, k)$ denote the corresponding cost benefit.
	\item Let $\tilde{x}$ denote the value of a temporary control policy. Let $\tilde{J}$ denote the corresponding cost benefit.
	\item Let $\lambda$ denote the coefficient of a quadratic regularization term, which is used to control the smoothness of the resulting control policy.
\end{itemize}

We are ready to present a dynamic programming algorithm for solving the CRP problem, which is detailed in the following CRP2 algorithm.

\begin{algorithm}[h] 
	\caption{CRP2}   
	\label{alg2}   
	{\bf Input:} an instance $\mathbb{M} = (A_0, I_0, T, \overline{x}, \mu, \delta_1, \delta_2, \alpha, \beta_1, \beta_2, \omega_1,$ $\omega_2)$ of the CRP model (5), three positive integers $N$, $M$, and $P$, regularization coefficient $\lambda$. \\
	{\bf Output:} a control policy $x$.	
	\begin{algorithmic}[1] 			
		\STATE // sentences 2-6 initialize all $x(i, j, k)$ and $J(i, j, k)$; // 
		\FOR {$i \leftarrow 0$ to $N$}
		\FOR {$j, k \leftarrow 0$ to $M$}
		\STATE $x(i, j, k) \leftarrow 0$; $J(i, j, k) \leftarrow 0$;
		\ENDFOR 
		\ENDFOR
		\STATE // sentences 8-23 calculate a complete optimal control policy; //
		\FOR {$i \leftarrow N-1$ down to $0$}
		\STATE // sentences 10-22 calculate an optimal control policy at time point $t_i$; //
		\FOR {$j,k \leftarrow 0$ to $M$}
		\STATE $\tilde{x} \leftarrow 0$; $\tilde{J} \leftarrow - \infty$;
		\FOR{$p \leftarrow 0$ to $P$}
		\STATE calculate $(A_{i+1}, I_{i+1})$ using the discretized version of Eqs. (3) with $(A_i, I_i)$ and $x_p$;
		\STATE  $j'\leftarrow \arg \underset{\tilde{j} \in \{0, 1, \cdots, M \}} \min |s_{\tilde{j}} - A_{i+1}|$; 
		\STATE  $k'\leftarrow \arg \underset{\tilde{k} \in \{0, 1, \cdots, M \}} \min |s_{\tilde{k}} - I_{i+1}|$; 
		\STATE $J \leftarrow \omega_2 I_{i+1} - \omega_1 x_p +\lambda x_p^2 + J(i+1, j', k') $;
		\IF {$J > \tilde{J}$}
		\STATE $\tilde{x} \leftarrow x_p$; $\tilde{J} \leftarrow J$;
		\ENDIF 		
		\ENDFOR
        \STATE $x(i,j,k) \leftarrow \tilde{x}$; $J(i,j,k) \leftarrow \tilde{J}$;
        \ENDFOR
        \ENDFOR
	\end{algorithmic} 
\end{algorithm} 

\begin{ex}
Let $\mathbb{M}_1$ denote the CRP instance given in Experiment 1. Again, let $x_1^*$ denote the CRP obtained in Experiment 1. Let $N = 50$, $M=2000$, $P = 100$, $\lambda = 0.1$. Running the CRP2 algorithm on $(\mathbb{M}_1, N, M, P, \lambda)$, we get a CRP, which is denoted $x_2^*$ and plotted in Fig. 5. 
\end{ex}

\begin{figure}[!ht]
	\setlength{\abovecaptionskip}{0.cm}
	\setlength{\belowcaptionskip}{-0.cm}
	\centering
	\includegraphics[scale=0.38]{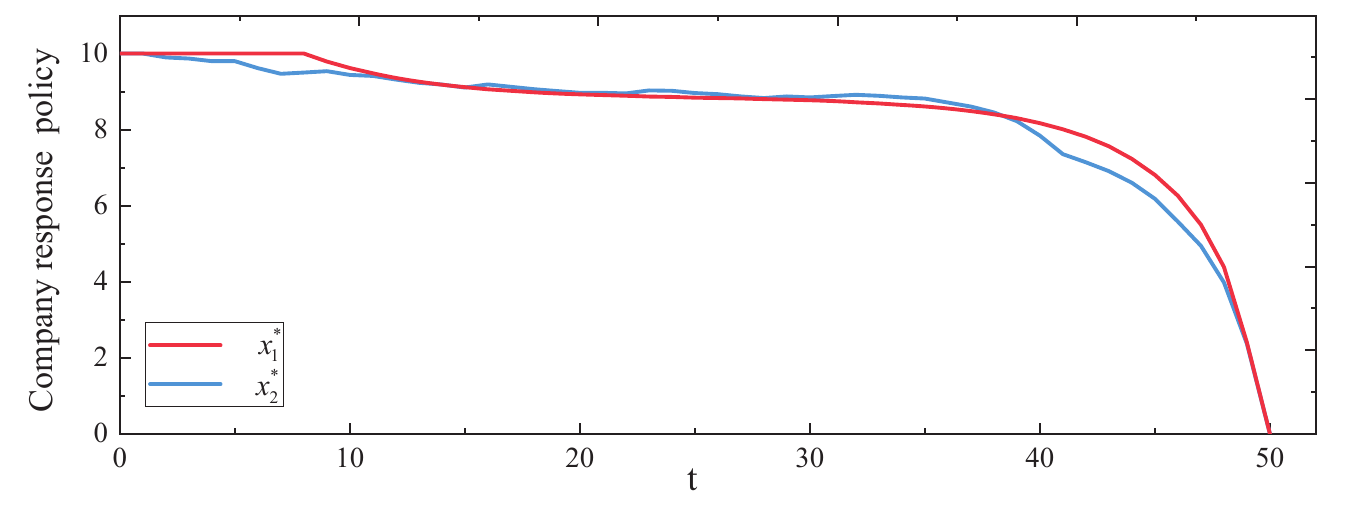}
	\caption{The CRP $x_2^*$ obtained in Experiment 4.}
\end{figure}

First, we observe that $J(x_2^*) = 4.02 \times 10^6$. In contrast, $J(x_1^*) = 4.06 \times 10^6$. As a result, $x_2^*$ is inferior to $J(x_1^*)$ in terms of the cost benefit. This is largely attributed to the rough discrete approximation of the CRP model. Second, the running time of the CRP algorithm in Experiment 1 is only less than 10 seconds. In contrast, the running time of the CRP2 algorithm in Experiment 4 is up to more than 6 hours. What is worse, with the increased refinement of discrete approximation of the CRP model, the time cost of the CRP2 algorithm would rise rapidly. Additionally, Fig. 5 shows that $x_2^*$ is much less smooth than $x_1^*$. Consequently, we conclude that the CRP2 algorithm is inferior to the CRP algorithm from three perspectives: the cost benefit of the resulting CRP, the time cost of running the algorithm, and the smoothness of the resulting CRP. Consequently, we recommend the CRP algorithm to serve as the basic algorithm for solving the CRP model.

\subsection{Implementation of the resulting company response policy}

The very key to implementing the recommended company response policy in real-world product co-creation activities lies in understanding the instance of the CRP model. It is seen from Eq. (6) that the instance involves 12 factors. Hence, the company needs to understand these factors. 

First, the initial community state can be observed by the company. Second, the product co-creation period and the maximal response rate are set directly by the company. Third, by collecting and averaging a collection of historical data on the community inflow rate and the two community outflow rates, the three rates can be estimated relatively accurately. Fourth, by collecting and averaging a collection of historical data on the inaction rate, the inaction rate may be estimated relatively accurately. Next, the standard cost and the standard benefit may be estimated relatively accurately by the company. Fifth, by collecting and statistically fitting a collection of historical data on the correspondence between the company response rate and the inactive-active conversion rate, the first influence functions may be approximated relatively accurately. Finally, by collecting and statistically fitting a collection of relevant historical data on the correspondence between the number of active participants and the inactive-active conversion rate, the second influence functions may be approximated relatively accurately. When all these tasks are accomplished, the recommended company response policy may be calculated and implemented.

\section{Further discussions}

In the preceding section, we validated the convergence and effectiveness of the CRP algorithm. In this section, we inspect the influence of some factors on the cost benefit of the CRP resulting from running the CRP algorithm (the resulting CRP, for short) through numerical experiments.

\subsection{Influence of the product co-creation period}

First, examine the influence of the product co-creation period on the cost benefit of the resulting CRP.

\begin{ex} 
Let $\mathcal{T} = \{100, 110, \cdots, 200\}$. Consider the following set of CRP instances.
\begin{equation*}
\begin{aligned}
  \mathbb{M}_T &= (100, 10000, T, 15, 12, 0.0001, 0.001, 0.1, \\
  & 0.05\arctan (0.3z), 0.01\ln (0.01z+1), 800, 20), T \in \mathcal{T}.
\end{aligned}
\end{equation*}
For each $T \in \mathcal{T}$, let $x_T^*$ denote the CRP resulting from running the CRP algorithm on $(\mathbb{M}_T, 10^{-6})$. Fig. 6(a) displays $J(x_T^*)$ versus $T$, $T \in \mathcal{T}$. We observe that $J(x_T^*)$ increases with the increase of $T$.
\end{ex}

We conduct totally $10^{4}$ similar experiments. In each of these experiments, we observe that the cost benefit of the resulting CRP increases with the increase of the product co-creation period. Hence, we conclude that, generally, the cost benefit of the resulting CRP is increasing with the product co-creation period. This conclusion informs that the extended product co-creation period implies the increased cost benefit of the co-creative product. Since the product co-creation period is typically preset by the community, it is not regarded as a key factor that helps enhance the cost benefit of product co-creation activity.

\begin{figure*}[!ht]
	\setlength{\abovecaptionskip}{0.cm}
	\setlength{\belowcaptionskip}{-0.cm}
	\centerline{\includegraphics[width=0.8\textwidth]{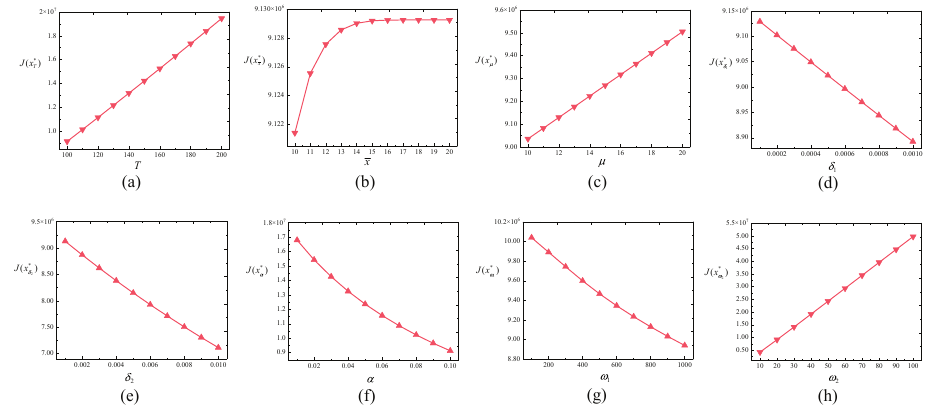}}
	\caption{The experimental results in (a) Experiment 5, (b) Experiment 6, (c) Experiment 7, (d) Experiment 8, (e) Experiment 9, (f) Experiment 10, (g) Experiment 11, and (h) Experiment 12.} 
\end{figure*}

\subsection{Influence of the maximal response rate}

Second, inspect the influence of the maximal response rate on the cost benefit of the resulting CRP.

\begin{ex} Let $\mathcal{X} = \{10, 11, \cdots, 20\}$. Consider the following set of CRP instances.
\begin{equation*}
\begin{aligned}
  \mathbb{M}_{\overline{x}} &= (100, 10000, 100, \overline{x}, 12, 0.0001, 0.001, 0.1, \\
  & 0.05\arctan (0.3z), 0.01\ln (0.01z+1), 800, 20), \overline{x} \in \mathcal{X}.
\end{aligned}
\end{equation*}
For each $\overline{x} \in \mathcal{X}$, let $x_{\overline{x}}^*$ denote the CRP resulting from running the CRP algorithm on $(\mathbb{M}_{\overline{x}}, 10^{-6})$. Fig. 6(b) exhibits $J(x_{\overline{x}}^*)$ versus $\overline{x}$, $\overline{x} \in \mathcal{X}$. We observe that, with the increase of $\overline{x}$, $J(x_{\overline{x}}^*)$ first increases rapidly and then flattens out. 
\end{ex}

We conduct totally $10000$ similar experiments. In each of these experiments, we observe that, with the increase of the maximal response rate, the cost benefit of the resulting CRP first increases rapidly and then flattens out. Hence, we conclude that this is a universal phenomenon. Below let us explain the phenomenon from two perspectives.

First, let $\overline{x}_1$ and $\overline{x}_2$ be a pair of maximal response rates with $\overline{x}_1 < \overline{x}_2$. It is easily seen from Eq. 1 that $X_{T, \overline{x}_1} \subset X_{T, \overline{x}_2}$. For any given $T$, $\mu$, $\delta_1$, $\delta_2$, $\alpha$, $\beta_1$, $\beta_2$, $\omega_1$, $\omega_2$, $A_0$, and $I_0$, let $x_1^{opt}$ and $x_2^{opt}$ be an optimal solution to the instance $\mathbb{M}_1 = (A_0, I_0, T, \overline{x}_1, \mu, \delta_1, \delta_2, \alpha, \beta_1, \beta_2, \omega_1, \omega_2)$ and the instance $\mathbb{M}_2 = (A_0, I_0, T, \overline{x}_2, \mu, \delta_1, \delta_2, \alpha, \beta_1, \beta_2, \omega_1, \omega_2)$, respectively. It is seen from the CRP model (5) that $J(x_1^{opt}) \leq J(x_2^{opt})$. Therefore, we can safely say the cost benefit of the resulting CRP is increasing with the increase of the maximal response rate.  

Second, it is seen from Figs. 2-4 that the response rate of the resulting CRP declines over time from the maximal response rate to zero. Furthermore, we observe that, with the increase of the maximal response rate, the declination begins at earlier time point. Hence, with the increase of the maximal response rate, the cost and, hence, the benefit of the resulting CRP tend to saturation. Therefore, with the increase of the maximal response rate, the cost benefit of the resulting CRP first increases rapidly and then flattens out. See Fig. 6(b).

In the physical world, there exists a maximum physically feasible response rate, which may be evaluated accurately by relevant experts. As a result, the actual maximal response rate is only physically achievable if it doesn't exceed the maximum physically feasible response rate. In the case where the restriction is met, the company may flexibly adjust the actual maximal response rate. In view of the asymptotic saturation of the cost effectiveness of the resulting CRP with the increase of the maximal response rate (see Fig. 6(b)), it is appropriate to choose a physically achievable maximal response rate at which the cost benefit of the resulting CRP flattens out.

\subsection{Influence of the community inflow rate}

Third, look into the influence of the community inflow rate on the cost benefit of the resulting CRP.

\begin{ex} Let $\mathcal{M} = \{10, 11, \cdots, 20\}$. Consider the following set of CRP instances.
\begin{equation*}
\begin{aligned}
  \mathbb{M}_{\mu} &= (100, 10000, 100, 15, \mu, 0.0001, 0.001, 0.1, \\
  & 0.05\arctan (0.3z), 0.01\ln (0.01z+1), 800, 20), \mu \in \mathcal{M}.
\end{aligned}
\end{equation*}
For each $\mu \in \mathcal{M}$, let $x_{\mu}^*$ denote the CRP resulting from running the CRP algorithm on $(\mathbb{M}_{\mu}, 10^{-6})$. Fig. 6(c) shows $J(x_{\mu}^*)$ versus $\mu$, $\mu \in \mathcal{M}$. We observe that $J(x_{\mu}^*)$ increases quickly with the increase of $\mu$. 
\end{ex}

We conduct totally $10000$ similar experiments. In each of these experiments, we observed that the cost benefit of the resulting CRP increases with the increase of the community inflow rate. Hence, we conclude that, generally, the cost benefit of the resulting CRP is increasing quickly with the increase of the community inflow rate. This conclusion informs that the lifted community inflow rate leads to the increased cost benefit of the co-creative product. 

Engaging as many potential customers as possible in online co-creation communities is the very key to realizing value co-creation. When it comes to product co-creation through company-sponsored online communities, more community participants implies more active participants, which in turn implies more valuable suggestions about the product design scheme, which in turn implies perfect product design scheme. Consequently, we suggest taking effective measures, such as investing more advertising expenditure on popular social media advertising platforms, to enlarge the customer base. Since the community inflow rate is relatively controllable and the resulting cost benefit is considerable, it is regarded as a key factor that helps enhance the cost benefit of product co-creation activity.

\subsection{Influence of the two community outflow rates}

Fourth, consider the influence of the two community outflow rates on the cost benefit of the resulting CRP.

\begin{ex} Let $\Delta_1 = \{0.0001, 0.0002, \cdots, 0.001\}$. Consider the following set of CRP instances.
\begin{equation*}
\begin{aligned}
  \mathbb{M}_{\delta_1} &= (100, 10000, 100, 15, 12, \delta_1, 0.001, 0.1, \\
  & 0.05\arctan (0.3z), 0.01\ln (0.01z+1), 800, 20), \delta_1 \in \Delta_1.
\end{aligned}
\end{equation*}
For each $\delta_1 \in \Delta_1$, let $x_{\delta_1}^*$ denote the CRP resulting from running the CRP algorithm on $(\mathbb{M}_{\delta_1}, 10^{-6})$. Fig. 6(d) showcases $J(x_{\delta_1}^*)$ versus $\delta_1$, $\delta_1 \in \Delta_1$. We observe that $J(x_{\delta_1}^*)$ decreases with the increase of $\delta_1$. 
\end{ex}

\begin{ex} Let $\Delta_2 = \{0.001, 0.002, \cdots, 0.01\}$. Consider the CRP instances
\begin{equation*}
\begin{aligned}
  \mathbb{M}_{\delta_2} &= (100, 10000, 100, 15, 12, 0.0001, \delta_2, 0.1, \\
  & 0.05\arctan (0.3z), 0.01\ln (0.01z+1), 800, 20), \delta_2 \in \Delta_2.
\end{aligned}
\end{equation*}
For each $\delta_2 \in \Delta_2$, let $x_{\delta_2}^*$ denote the CRP resulting from running the CRP algorithm on $(\mathbb{M}_{\delta_2}, 10^{-6})$. Fig. 6(e) showcases $J(x_{\delta_2}^*)$ versus $\delta_2$, $\delta_2 \in \Delta_2$. We observe that $J(x_{\delta_2}^*)$ decreases with the increase of $\delta_2$. 
\end{ex}

We conduct totally $10000$ similar experiments. In each of these experiments, we observe that the cost benefit of the resulting CRP decreases with the increase of either of the two community outflow rates. Hence, we conclude that, generally, the cost benefit of the resulting CRP is decreasing with the increase of either of the two community outflow rates. This conclusion informs that the two reduced community outflow rates imply the increased cost benefit of the co-creative product. 

Indubitably, the massive outflow of participants from a company-sponsored online community is a heavy blow to value co-creation. After all, maintaining the stability of existing participants contributes to the emergence of numerous valuable suggestions, leading to high-quality product design scheme. Consequently, we suggest taking effective measures, such as paying for valuable suggestions, to detain the community participants. Since the two outflow rates are controllable, they are regarded as key factors that help enhance the cost benefit of product co-creation activity.

\subsection{Influence of the inaction rate}

Fifth, consider the influence of the inaction rate on the cost benefit of the resulting CRP.

\begin{ex} Let $\mathcal{A} = \{0.01, 0.02, \cdots, 0.1\}$. Consider the following set of CRP instances.
\begin{equation*}
\begin{aligned}
  \mathbb{M}_{\alpha} &= (100, 10000, 100, 15, 12, 0.0001, 0.001, \alpha, \\
  & 0.05\arctan (0.3z), 0.01\ln (0.01z+1), 800, 20), \in \mathcal{A}.
\end{aligned}
\end{equation*}
For each $\alpha \in \mathcal{A}$, let $x_{\alpha}^*$ denote the CRP resulting from running the CRP algorithm on $(\mathbb{M}_{\alpha}, 10^{-6})$. Fig. 6(f) depicts $J(x_{\alpha}^*)$ versus $\alpha$, $\alpha \in \mathcal{A}$. We observe that $J(x_{\alpha}^*)$ decreases with the increase of $\alpha$. 
\end{ex}

We conduct totally $10000$ similar experiments. In each of these experiments, we observe that the cost benefit of the resulting CRP decreases with the increase of the inaction rate. Hence, we conclude that, generally, the cost benefit of the resulting CRP is decreasing with the increase of the inaction rate. This conclusion informs that the reduced inaction rate leads to the increased benefit of product co-creation. 

Around 400 years BC, Sophocles, the famous Ancient Greek playwright, wrote the aphorism: Wisdom outweighs any wealth. As the combination of experience, knowledge, and sensible decision, wisdom is the key to wealth. With regard to product co-creation through company-sponsored online community, taking full advantage of the participants' wisdom undoubtedly helps to yield perfect product design scheme and hence great wealth. Consequently, we suggest taking effective measures, such as paying active participants for their excellent ideas, to unearth their potential.

\subsection{Influence of the standard cost}

Next, consider the influence of the standard cost on the cost benefit of the resulting CRP.

\begin{ex} Let $\mathcal{W}_1 = \{100, 200, \cdots, 1000\}$. Consider the following set of CRP instances.
\begin{equation*}
\begin{aligned}
  \mathbb{M}_{\omega_1} &= (100, 10000, 100, 15, 12, 0.0001, 0.001, 0.1, \\
  & 0.05\arctan (0.3z), 0.01\ln (0.01z+1), \omega_1, 20), \omega_1 \in \mathcal{W}_1.
\end{aligned}
\end{equation*}
For each $\omega_1 \in \mathcal{W}_1$, let $x_{\omega_1}^*$ denote the CRP resulting from running the CRP algorithm on $(\mathbb{M}_{\omega_1}, 10^{-6})$. Fig. 6(g) draws $J(x_{\omega_1}^*)$ versus $\omega_1$, $\omega_1 \in \mathcal{W}_1$. We observe that $J(x_{\omega_1}^*)$ decreases with the increase of $\omega_1$. 
\end{ex}

We conduct totally $10000$ similar experiments. In each of these experiments, we observe that the cost benefit of the resulting CRP decreases with the increase of the standard cost. Hence, we conclude that, generally, the cost benefit of the resulting CRP is decreasing with the increase of the standard cost. This conclusion informs that the reduced standard cost implies the increased cost benefit of the co-creative product. 

With everything in business, the benefits gained should exceed the cost incurred. Product co-creation is no exception. While keeping the high intensity of interacting with the community participants, we suggest taking necessary measures, such as improving the operational efficiency of the relevant department, to reduce the cost for company response.

\subsection{Influence of the standard benefit}

Finally, examine the influence of standard benefit on the cost benefit of the resulting CRP.

\begin{ex} Let $\mathcal{W}_2 = \{10, 20, \cdots, 100\}$. Consider the following set of CRP instances.
\begin{equation*}
\begin{aligned}
  \mathbb{M}_{\omega_2} &= (100, 10000, 100, 15, 12, 0.0001, 0.001, 0.1, \\
  & 0.05\arctan (0.3z), 0.01\ln (0.01z+1), 800, \omega_2), \omega_2 \in \mathcal{W}_2.
\end{aligned}
\end{equation*}
For each $\omega_2 \in \mathcal{W}_2$, let $x_{\omega_2}^*$ denote the CRP resulting from running the CRP algorithm on $(\mathbb{M}_{\omega_2}, 10^{-6})$. Fig. 6(h) depicts $J(x_{\omega_2}^*)$ versus $\omega_2$, $\omega_2 \in \mathcal{W}_2$. We observe that $J(x_{\omega_2}^*)$ increases with the increase of $\omega_2$. 
\end{ex}

We conduct totally $10000$ similar experiments. In each of these experiments, we observe that the cost benefit of the resulting CRP increases with the increase of the standard benefit. Hence, we conclude that, generally, the cost benefit of the resulting CRP is increasing with the increase of the standard benefit. This conclusion informs that the lifted standard benefit leads to the increased cost benefit of the co-creative product. 

One disruptive innovation outvalues one thousand minor betterments. Generous rewards rouse one to heroism. The recipe for increasing the standard benefit is to simulate active participants' potential to the maximum extent by posting a high reward. Consequently, we suggest awarding active participants differentially based on the quality of their ideas, with the intention of arousing their passion for pursuing top-level product design scheme.

\begin{rk}
It is of practical importance to ask if the performance of the CRP algorithm is continuously dependent on the model parameters. It can be seen from the above discussions that the answer is yes. Therefore, it is expected that, with the increasingly accurate observations of these parameters, the performance of the CRP algorithm would become increasingly predictable.
\end{rk}

\section{Concluding remarks}

In the context of product co-creation through company-sponsored online community, the problem of finding a cost-effective company response policy has been proposed. This problem has been reduced to an optimal control problem. An algorithm for solving the latter problem has been presented. The convergence and effectiveness of the algorithm has been corroborated. Consequently, the company response policy resulting from running the algorithm have been recommended.

There are several relevant problems to be resolved. First, gathering numerous realistic CRP model-related data helps to yield cases of calculating and implementing the recommended CRP. Second, early contributions of a participant may inspire others to become intrigued and look to make additional contributions. In the future modification of the community state evolutionary model, this delayed influence should be incorporated. Thirdly, since product co-creation based on company-sponsored online community may be viewed as a cooperative game between the company and the community participants, it is worthwhile to study the original problem from the perspective of cooperative game theory \cite{Basar2012}. Next, it is worth exploring the formation mechanism of an autonomous online co-creation community through epidemic modeling \cite{Zino2022}. Finally, the methodology developed in this paper may be borrowed to address some other issues, such as cost-effective advertising \cite{HuangKF2022} and cost-effective cyber defense \cite{YangLX2021a, YangLX2021b}.


%



\ifCLASSOPTIONcaptionsoff
  \newpage
\fi




\begin{thebibliography}{67}

\bibitem[1]{Prahalad2004a}
C. K. Prahalad and V. Ramaswamy, ``Co-creating unique value with customers'', \emph{Strategy \& Leadership}, vol. 32, no. 3, pp. 4-9, 2004.

\bibitem[2]{Prahalad2004b}
C. K. Prahalad and V. Ramaswamy, \emph{The Future of Competition: Co-Creating Unique Value with Customers}, Harvard Business Review Press, 2004.

\bibitem[3]{Ramaswamy2010}
V. Ramaswamy and F. Gouillart, \emph{The Power of Co-Creation: Build It with Them to Boost Growth, Productivity, and Profits}, Free Press, 2010.

\bibitem[4]{Ramaswamy2014}
V. Ramaswamy and K. Ozcan, \emph{The Co-Creation Paradigm}, Stanford Business Books, 2014.

\bibitem[5]{Faraj2011}
S. Faraj, S. L. Jarvenpaa, and A. Majchrzak, ``Knowledge collaboration in online communities'', \emph{Organization Science}, vol. 22, no. 5, pp. 1224-1239, 2011.

\bibitem[6]{Lusch2015}
R. F. Lusch and S. Nambisan, "Service innovation: A service-dominant logic perspective", \emph{MIS Quarterly}, Vol. 39, no. 1, pp. 155-176, 2015.

\bibitem[7]{Nambisan2017}
S. Nambisan, K. Lyytinen, A. Majchrzak, and M. Song, ``Digital innovation management:
Reinventing innovation management research in a digital world'', MIS Quarterly, vol. 41, no.
1, pp. 223-238, 2017.

\bibitem[8]{YanJ2018}
J. Yan, D. Leidner, and H. Benbya, ``The differential Innovativeness outcomes of user and employee participation in an online user innovation community'', \emph{Journal of Management Information Systems}, vol. 35, no. 3, pp. 1–34, 2018.

\bibitem[9]{Gangi2009}
P. M. Di Gangi and M. Wasko, ``Organizational adoption of user innovation from a user innovative community: A case study of Dell IdeaStorm'', \emph{Decision Support Systems}, vol. 48, no. 1, pp. 303-312, 2009.

\bibitem[10]{Wong2016}
T. Y. Wong, G. Peko, D. Sundaram, and S. Piramuthu, ``Mobile environments and innovative co-reation processes \& ecosystems'', \emph{Information \& Management}, vol. 53, no. 3, pp. 336-344, 2016.

\bibitem[11]{Tavakoli2017}
A. Tavakoli, D. Schlagwein, and D. Schoder, ``Open strategy: Literature review, re-analysis of cases and conceptualisation as a practice'', \emph{The Journal of Strategic Information Systems}, vol. 26, no. 3, pp. 163-184, 2017.

\bibitem[12]{Essamri2019}
A. Essamri, S. MeKechnie, and H. Winklhofer, ``Co-creating corporate brand identity with online brand communities: A managerial perspective'', \emph{Journal of Business Research}, vol. 96, pp. 366-375, 2019.

\bibitem[13]{Bhalla2011}
G. Bhalla, \emph{Collaboration and Co-creation: New Platforms for Marketing and Innovation}, Springer Science+Business Media, 2011.

\bibitem[14]{Zwass2010}
V. Zwass, ``Co-creation: Toward a taxonomy and an integrated research perspective'', \emph{International Journal of Electronic Commerce}, vol. 15, no. 1, pp. 11-48, 2010.

\bibitem[15]{Porter2008}
C. E. Porter and N. Donthu, ``Cultivating trust and harvesting value in virtual communities'', \emph{Management Science}, vol. 54, no. 1, pp. 113-128, 2008.

\bibitem[16]{Faraj2015}
S. Faraj, S. Kudaravalli, and M. Wasko, ``Leading collaboration in online communities'', \emph{MIS Quarterly}, vol. 32, no. 2, pp. 393-412, 2015.

\bibitem[17]{Constantinides2015}
E. Constantinides, L. A. Brunink, and C. Lorenzo-Romero, ``Customer motives and benefits for participating in online co-creation activities'', \emph{International Journal of Internet Marketing and Advertising}, vol. 9, no. 1, pp. 21-48, 2015.

\bibitem[18]{Fernandes2016}
T. Fernandes and P. Remelhe, ``How to engage customers in co-creation: customers' motivations for collaborative innovation'', \emph{Journal of Strategic Marketing}, vol. 29, no. 1, pp. 218-244, 2019.

\bibitem[19]{Spaulding2010}
T. J. Spaulding, ``How can virtual communities create value for business'', \emph{Management Science}, vol. 54, no. 1, pp. 113-128, 2008.

\bibitem[20]{Fisher2019}
G. Fisher, ``Online communities and firm advantages'', \emph{Academy of Management Review}, vol. 44, no. 2, pp. 279-298, 2019.

\bibitem[21]{Rubio2020}
N. Rubio, N. Villaenor, and M. J. Yague, ``Sustainable co-creation behavior in a virtual community: Antecedents and moderating effect of participants's perception of own expertise'', \emph{Sustainability}, vol. 12, art. no. 8151, 2020.

\bibitem[22]{Jeppesen2006}
L. B. Jeppesen and L. Frederiksen, ``Why do users contribute to firm-hosted user communities? The case of computer-controlled music instruments'', \emph{Organization Science}, vol. 17, no. 1, pp. 45-63, 2006. 

\bibitem[23]{Franke2013}
N. Franke, P. Keinz, and K. Klausberger, ``Does this sound like a fair deal?'': Antecedents and consequences of fairness expectations in the individual's decision to participate in firm innovation'', \emph{Organization Science}, vol. 24, no. 5, pp. 1495-1516, 2013.
 
\bibitem[24]{Manchanda2015}
P. Manchanda, G. Packard, and A. Pattabhiramaiah, ``Social dollars: The economic impact of customer participation in a firm-sponsored online customer community'', \emph{Marketing Science}, vol. 34, no. 3, pp. 309-472, 2015.

\bibitem[25]{Priharsari2020}
D. Priharsari, B. Abedin, and E. Mastio,  ``Value co-creation in firm sponsored online communities: What enables, constraints, and shapes value'', \emph{Internet Research}, vol. 30, no. 3, pp. 763-788, 2020.

\bibitem[26]{ChenL2014}
L. Chen, J. R. Marsden, and Z. Zhang, ``Theory and analysis of company-sponsored value co-creation'', \emph{Journal of Management Information Systems}, vol. 34, no. 3, pp. 141-172, 2014.

\bibitem[27]{Dietz2002}
K. Dietz and J. A. P. Heesterbeek, ``Daniel Bernoulli's epidemiological model revisited'', \emph{Mathematical Biosciences}, vol. 180, pp. 1-21, 2002.

\bibitem[28]{Kermack1927}
W. O. Kermack and A. G. McKendrick, ``A contribution to the mathematical theory of epidemics'',
\emph{Proceedings of the Royal Society of London, Series A, Containing Papers of a Mathematical and Physical Character}, vol. 115, no. 772, pp. 700-721, 1927.

\bibitem[29]{Kermack1933}
W. O. Kermack and A. G. McKendrick, ``A contribution to the mathematical theory of epidemics III: Further studies of the problem of endemicity'', \emph{Proceedings of the Royal Society A: Mathematical, Physical, and Engineering Sciences}, vol. 141, no. 83, pp. 94-122, 1933.

\bibitem[30]{Keeling2005}
M. J. Keeling and K. T. D. Eames, ``Networks and epidemic models'', \emph{Journal of Royal Society Interface}, vol. 2, no. 4, pp. 295-307, 2005.

\bibitem[31]{Pare2020}
P. E. Pare, C. L. Beck, and T. Basar, ``Modeling, estimation, and analysis of epidemics over networks: An overview'', \emph{Annual Reviews in Control}, vol. 50, pp. 345-360, 2020.

\bibitem[32]{Allen2008}
L. J. S. Allen, ``An introduction to stochastic epidemic models'', in: F. Brauer, P. van den Driessch, and J. Wu (eds), \emph{Mathematical Epidemiology}, Lecture Notes in Mathematics, vol 1945, Springer, Berlin, Heidelberg, 2008.

\bibitem[33]{Britton2010}
T. Britton, ``Stochastic epidemic models: a survey'', \emph{Mathematical Biosciences}, vol. 225, no. 1, pp. 24-35, 2010.

\bibitem[34]{Nowzari2017}
C. Nowzari, V. M. Preciado, and G. J. Pappas, ``Optimal resource allocation for control
of networked epidemic models'', \emph{IEEE Transactions on Control of Network Systems},
vol. 4, no. 2, pp. 159-169, 2017.

\bibitem[35]{Somers2021}
V. L. J. Somers and I. R. Manchester, ``Resource allocation for control of
spreading processes via convex optimization'', \emph{IEEE Control Systems Letters}, vol. 5,
no. 2, pp. 547-552, 2021.

\bibitem[36]{Pare2018}
P. E. Pare, C. L. Beck, and A. Nedic, ``Processes over time-varying networks'', \emph{IEEE Transactions on Control of Network Systems}, vol. 5, no. 3, pp. 1322-1334, 2018.

\bibitem[37]{Leitch2019}
J. Leitch, K. A. Alexander, and S. Sengupta, ``Epidemic thresholds on temporal
networks: a review and open questions'', \emph{Applied Network Science}, vol. 4, no. 1, art. no. 105, 2019.

\bibitem[38]{Zino2022}
L. Zino and M. Cao, ``Analysis, prediction, and control of epidemics: A survey from scalar to dynamic network models'', \emph{IEEE Circuits and Systems Magazine}, vol. 21, no. 4, pp. 4-23, 2021.

\bibitem[39]{Hsieh2004}
Y.-H. Hsieh, J.-Y. Lee, and H.-L. Chang, ``SARS epidemiology modeling'', \emph{Emerging Infectious Diseases}, vol. 10, no. 6, pp. 1165-1167, 2004.

\bibitem[40]{Gatto2020}
M. Gatto, E. Bertuzzo, L. Mari, S. Miccoli, L. Carraro, R. Casagrandi, and A. Rinaldo,
``Spread and dynamics of the COVID-19 epidemic in Italy: Effects of emergency containment measures'', \emph{Proceedings of the National Academy of Sciences}, vol. 117, no. 19,
pp. 10484-10491, 2020.

\bibitem[41]{JiangY2015}
Y. Jiang and J. C. Jiang, ``Diffuson in social networks: a multiagent perspective'', \emph{IEEE Transactions on Systems, Man, and Cybernetics: Systems}, vol. 45, no. 2, pp. 198-213, 2015.

\bibitem[42]{XingL2013}
L. Xing, G. Levitin, C. Wang, and Y. Dai, ``Reliability of systems subject to failures with dependent propagation effect'', \emph{IEEE Transactions on Systems, Man, and Cybernetics: Systems}, vol. 43, no. 2, pp. 277-290, 2013.

\bibitem[43]{Goffman1964}
W. Goffman and V. A. Newill, ``Generalization of epidemic Theory: An application to the transmission of ideas'', \emph{Nature}, vol. 204, pp. 225-228, 1964.

\bibitem[44]{Daley1964}
D. J. Daley and D. G. Kendall, ``Epidemics and rumours'', \emph{Nature}, vol. 204, p. 1118, 1964.

\bibitem[45]{YangLX2018a}
L. X. Yang, P. Li, X. Yang, Y. Wu, and Y. Y. Tang, ``On the competition of two conflicting messages'', \emph{Nonlinear Dynamics}, vol. 91, no. 3, pp. 1853-1869, 2018.

\bibitem[46]{Kephart1991}
J. O. Kephart and S. R. White, ``Directed-graph epidemiological models of computer viruses'', in: \emph{Proceedings of 1991 IEEE  Computer Society Symposium on Research in Security and Privacy}, pp. 343-359, 1991.

\bibitem[47]{YuS2015}
S. Yu, G. Gu, A. Barnawi, S. Guo, and I. Stojmenovic, ``Malware propagation in large-scale networks'', \emph{IEEE Transactions on Knowledge and Data Engineering}, vol. 27, no. 1, pp. 170-179, 2015.

\bibitem[48]{LiP2018a}
P. Li, X. Yang, L. X. Yang, Q. Xiong, Y. Wu, and Y. Y. Tang, ``The modeling and analysis of the word-of-mouth marketing'', \emph{Physica A: Statistical Mechanics and its Applications}, vol. 493, pp. 1-16, 2018.

\bibitem[49]{ZhangTR2019}
T. Zhang, P. Li, L. X. Yang, X. Yang, Y. Y. Tang, and Y. Wu, ``A discount strategy in word-of-mouth marketing'', \emph{Communications in Nonlinear Science and Numerical Simulation}, vol. 74, pp. 167-179, 2019.

\bibitem[50]{Liberzon2012}
D. Liberzon, \emph{Calculus of Variations and Optimal Control Theory: A Concise Introduction}, Princeton University Press, 2012.

\bibitem[51]{Aldila2020}
D. Aldila, M. Z. Ndii, and B. M. Samiadji, ``Optimal control on COVID-19 eradication program in Indonesia under the effect of community awareness'', \emph{Mathematical Biosciences and Engineering}, vol. 17, no. 6, pp. 6355-6389, 2020. 

\bibitem[52]{Sameni2022}
R. Sameni, ``Model-based prediction and optimal control of pandemics by non-pharmaceutical interventions'', \emph{IEEE Journal of Selected Topics in Signal Processing}, DOI: 10.1109/JSTSP.2021.3129118. 

\bibitem[53]{Eshghi2016}
S. Eshghi, M. H. R. Khouzani, S. Sarkar, and S. S. Venkatesh, ``Optimal patching in clustered malware epidemics'', \emph{IEEE/ACM Transactions on Networking}, vol. 24, no. 1, pp. 283-298, 2016.

\bibitem[54]{YangLX2021a}
L. X. Yang, P. Li, X. Yang, Y. Xiang, F. Jiang, and W. Zhou, ``Effective quarantine and recovery scheme against advanced persistent threat'', \emph{IEEE Transactions on Systems, Man, and Cybernetics: Systems}, vol. 51, no. 10, pp. 5977-5991, 2021.

\bibitem[55]{ZhaoJ2019}
J. Zhao, L. X. Yang, X. Zhong, X. Yang, Y. Wu, and Y. Y. Tang, ``Minimizing the impact of a rumor via isolation and conversion'', \emph{Physica A: Statistical Mechanics and its Applications}, vol. 526, art. no. 120867, 2019.

\bibitem[56]{LinY2021}
Y. Lin, X. Wang, F. Hao, Y. Jiang, Y. Wu, G. Min, D. He, S. Zhu, and W. Zhao, ``Dynamic control of fraud information spreading in mobile social networks,'' \emph{IEEE Transactions on Systems, Man, and Cybernetics: Systems}, vol. 51, no. 6, 2021, pp. 3725-3738, 2021.

\bibitem[57]{Basar2012}
T. Basar, \emph{Dynamic Noncooperative Game Theory}, Academic Press, 2012.

\bibitem[58]{HuangDW2021}
D. W. Huang, L. X. Yang, P. Li, X. Yang, and Y. Y. Tang, ``Developing cost-effective rumor-refuting strategy through game-theoretic approach,'' \emph{IEEE Systems Journal}, vol. 15, no. 4, pp. 5034-5045, 2021.

\bibitem[59]{YangLX2019a} 
L. X. Yang, P. Li, Y. Zhang, X. Yang, Y. Xiang, and Y. Y. Tang, ``Effective repair strategy against advanced persistent threat: A differential game approach'', \emph{IEEE Transactions on Information Forensics and Security}, vol. 14, no. 7, pp. 1713-1728, 2019.

\bibitem[60]{YangLX2021b}
L. X. Yang, K. Huang, X. Yang, Y. Zhang, Y. Xiang, Y. Y. Tang, ``Defense against advanced persistent threat through data backup and recovery'', \emph{IEEE Transactions on Network Science and Engineering}, vol. 8, no. 3, pp. 2001-2013, 2021.

\bibitem[61]{McAsey2012}
M. McAsey, L. Mou, and W. Han, ``Convergence of the forward-backward sweep method in optimal control'', \emph{Computational Optimization and Applications}, vol. 53, pp. 207-226, 2012.

\bibitem[62]{Sethi2000}
S. P. Sethi and G. L. Thompson, \emph{Optimal Control Theory: Applications to Management Science and Economics (2nd edition)}, Springer, New York, 2000.

\bibitem[63]{Burden2015}
R. L. Burden, J. D. Faires, and A. M. Burden, \emph{Numerical Analysis, 10th Edition}, Cengage Learning, 2015.

\bibitem[64]{Bertsekas2012}
D. Bertsekas, \emph{Dynamic Programming and Optimal Control}, 4th Edition, Athena Scientific, 2012.

\bibitem[65]{HuangKF2022}
K. Huang, L. X. Yang, X. Yang, and Y. Y. Tang, ``Effective multiplatform advertising policy'', \emph{IEEE Transactions on Systems, Man, and Cybernetics: Systems}, vol. 52, no. 7, pp. 4483-4493, 2022.

\end{thebibliography}
%

\begin{IEEEbiography}
[{\includegraphics[width=1in,height=1.25in,clip,keepaspectratio]{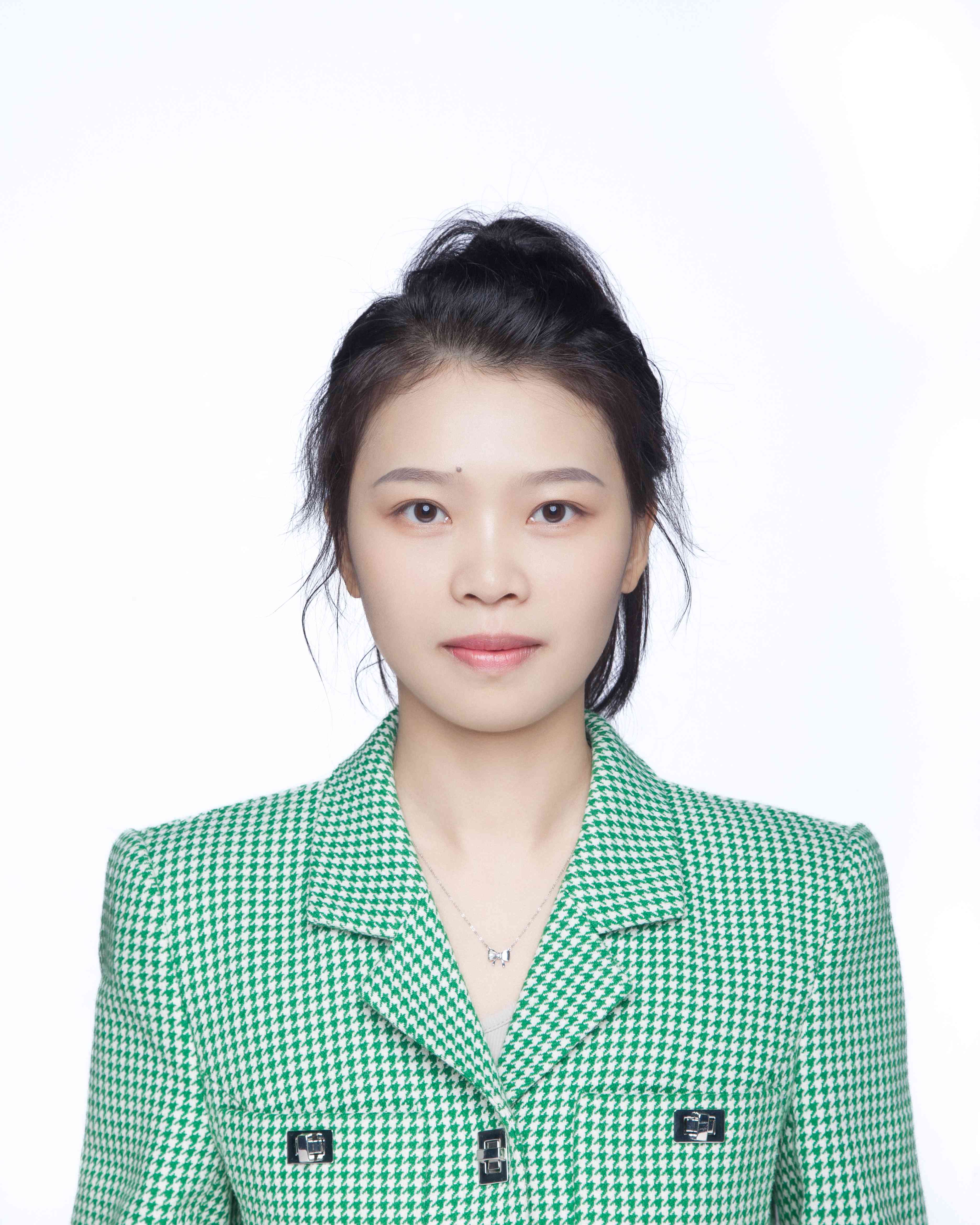}}]{Jiamin Hu} received the B.E. degree in Internet of Things Engineering from Chongqing University of Technology, Chongqing, China, in 2017. She is currently pursuing the Ph.D. degree in Chongqing University. She has published two academic papers in peer-reviewed international journals. Her research interests include applied optimal control theory and cyber security.
\end{IEEEbiography}

\begin{IEEEbiography}
[{\includegraphics[width=1in,height=1.25in,clip,keepaspectratio]{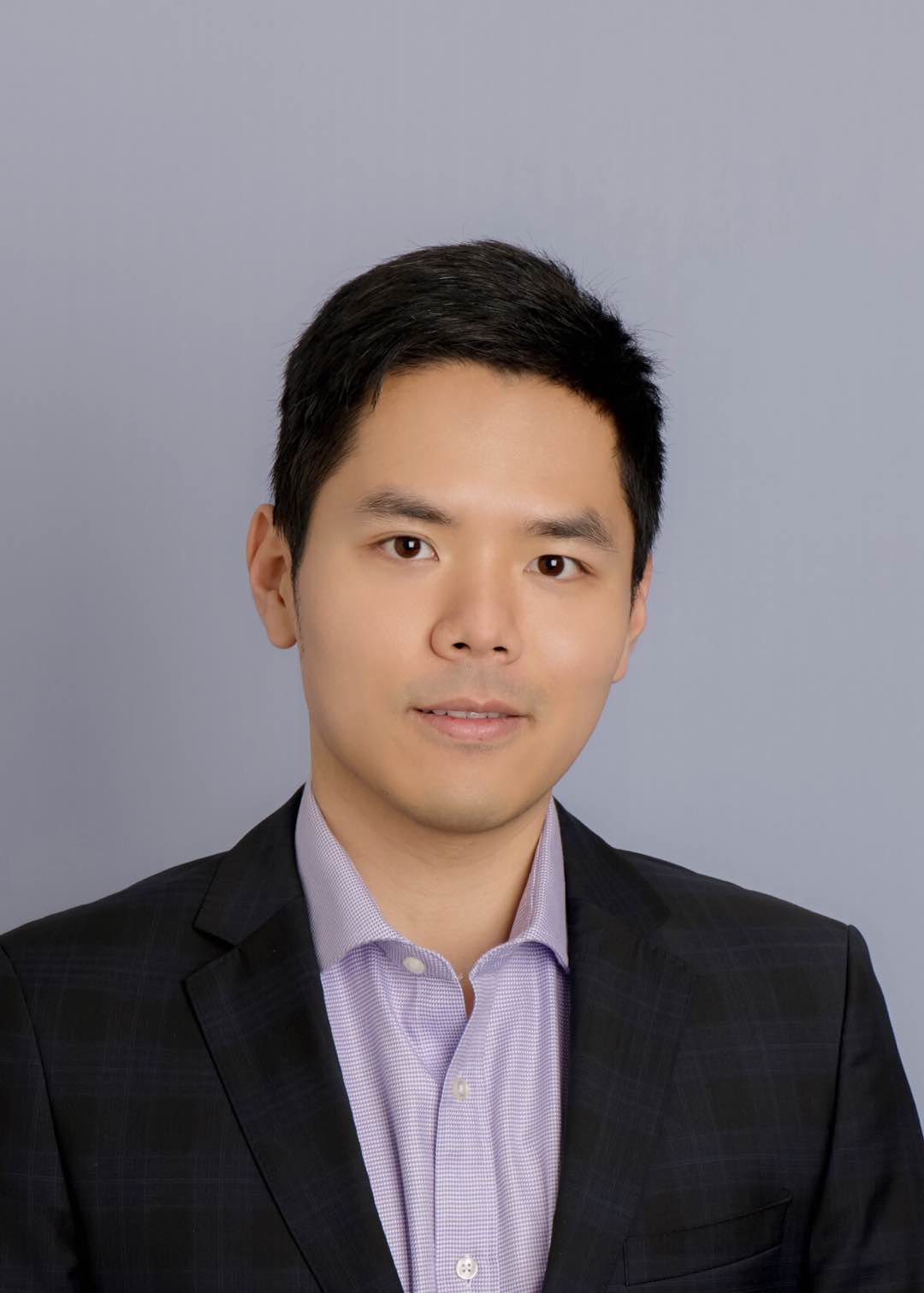}}]{Lu-Xing Yang} (Member, IEEE) received the B.Sc. degree in mathematics and the Ph.D. degree in computer science from Chongqing University, in 2012 and 2015, respectively. He is a Lecturer with the School of Information Technology, Deakin University. He has visited Imperial College London, from 2014 to 2015. He was a Post-Doctoral Research Fellow with the Delft University of Technology, from 2016 to 2017. He was an Alfred Post-Doctoral Research Fellow with Deakin University, from 2018 to 2020. He has published more than 60 academic articles in peer-reviewed international journals. His research interests include epidemic dynamics and cybersecurity.
\end{IEEEbiography}

\begin{IEEEbiography}
[{\includegraphics[width=1in,height=1.25in,clip,keepaspectratio]{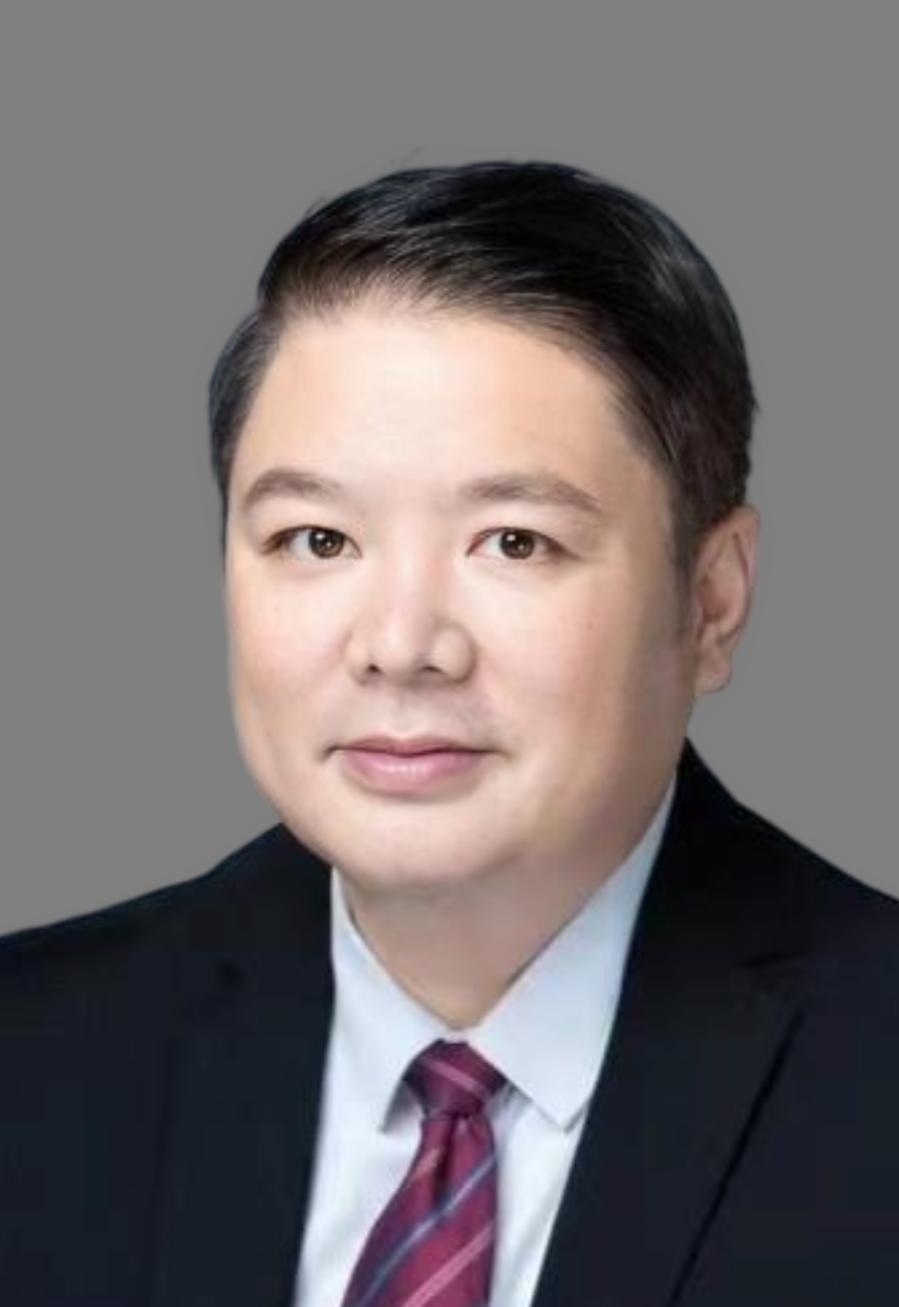}}]{Xiaofan Yang} (Member, IEEE) is a Professor with the School of Big Data and Software Engineering, Chongqing University. He received the B.Sc. degree in mathematics from Sichuan University in 1985, and the M.Sc. degree in mathematics and the Ph.D. degree in computer science from Chongqing University, in 1988 and 1994, respectively. He He has visited the University of Reading, from 1998 to 1999. He has published more than 160 academic articles in peer-reviewed journals. More than 20 students have received the Ph.D. degree under his supervision. His research interests include fault-tolerant computing, epidemic dynamics, and cybersecurity.
\end{IEEEbiography}

\begin{IEEEbiography}
[{\includegraphics[width=1in,height=1.25in,clip,keepaspectratio]{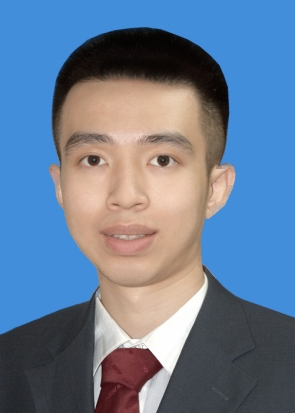}}]{Kaifan Huang} received the B.E. in software engineering from in 2018 in School of Big Data and Software Engineering from Chongqing University, where he is currently working toward the Ph.D. degree. He has published eight academic papers in peer-reviewed international journals. His research interests include viral marketing, rumor propagation, and cyber security.
\end{IEEEbiography}

\begin{IEEEbiography}
[{\includegraphics[width=1in,height=1.25in,clip,keepaspectratio]{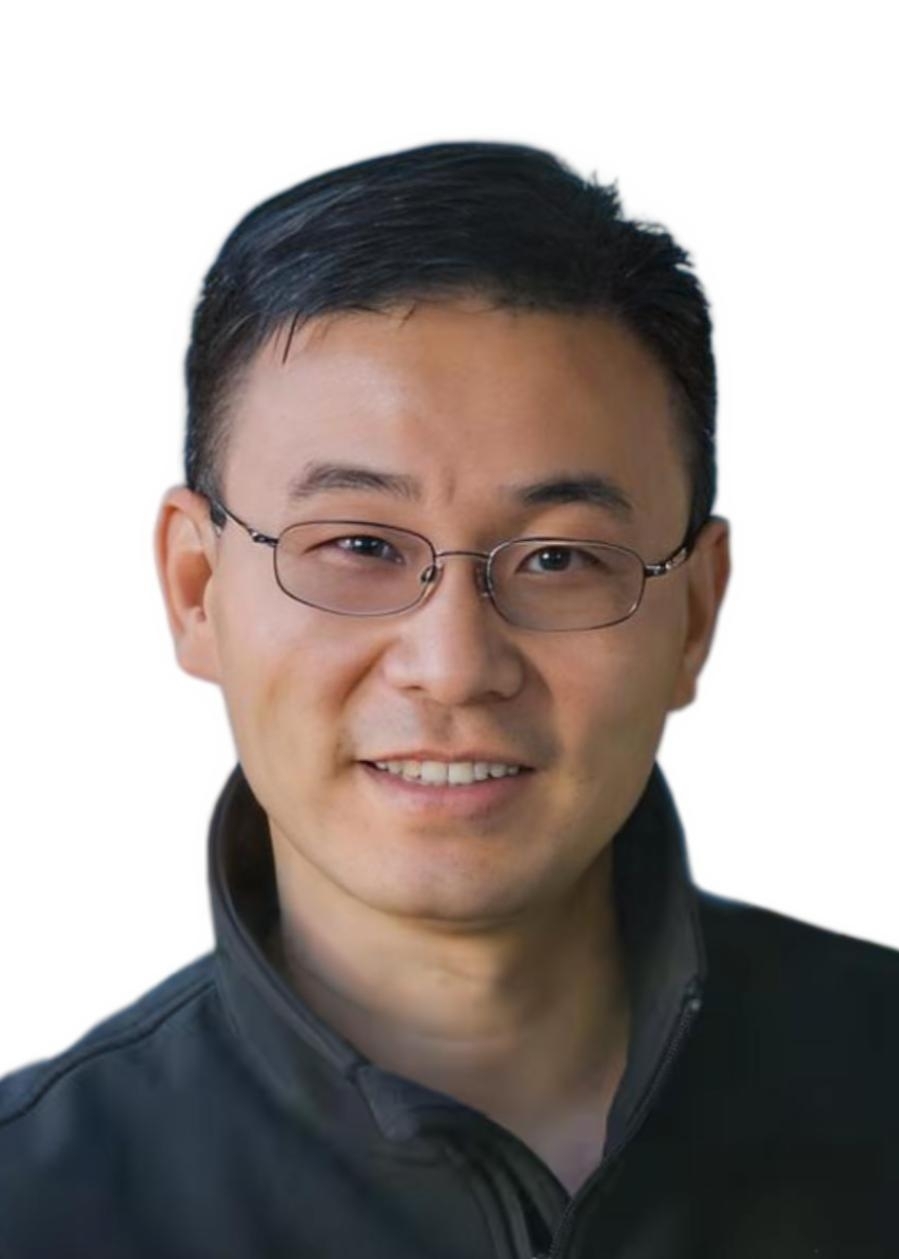}}]{Gang Li} (Senior Member, IEEE) is a Professor with the School of Information Technology, Deakin University, Australia. His research interests include data mining, data privacy, causal discovery, and business intelligence. He received the B.E. and M.Eng. degrees in computer science from the Xi’an Shiyou University, Xi’an, Shaanxi, China. He received the PH.D. degree in  computer science from the Deakin University. Dr. Li serves as the Chair for IEEE Task Force on Educational Data Mining, from 2020 to 2023, and the IEEE Enterprise Information Systems Technical Committee; and the Vice Chair for the IEEE Data Mining and Big Data Analytics Technical Committee, from 2017 to 2019. He is an Associate Editor of various A/A* journals, including Decision Support Systems (Elsevier), Journal of Travel Research (Sage), and Cybersecurity (Springer).
\end{IEEEbiography}

\begin{IEEEbiography}
[{\includegraphics[width=1in,height=1.25in,clip,keepaspectratio]{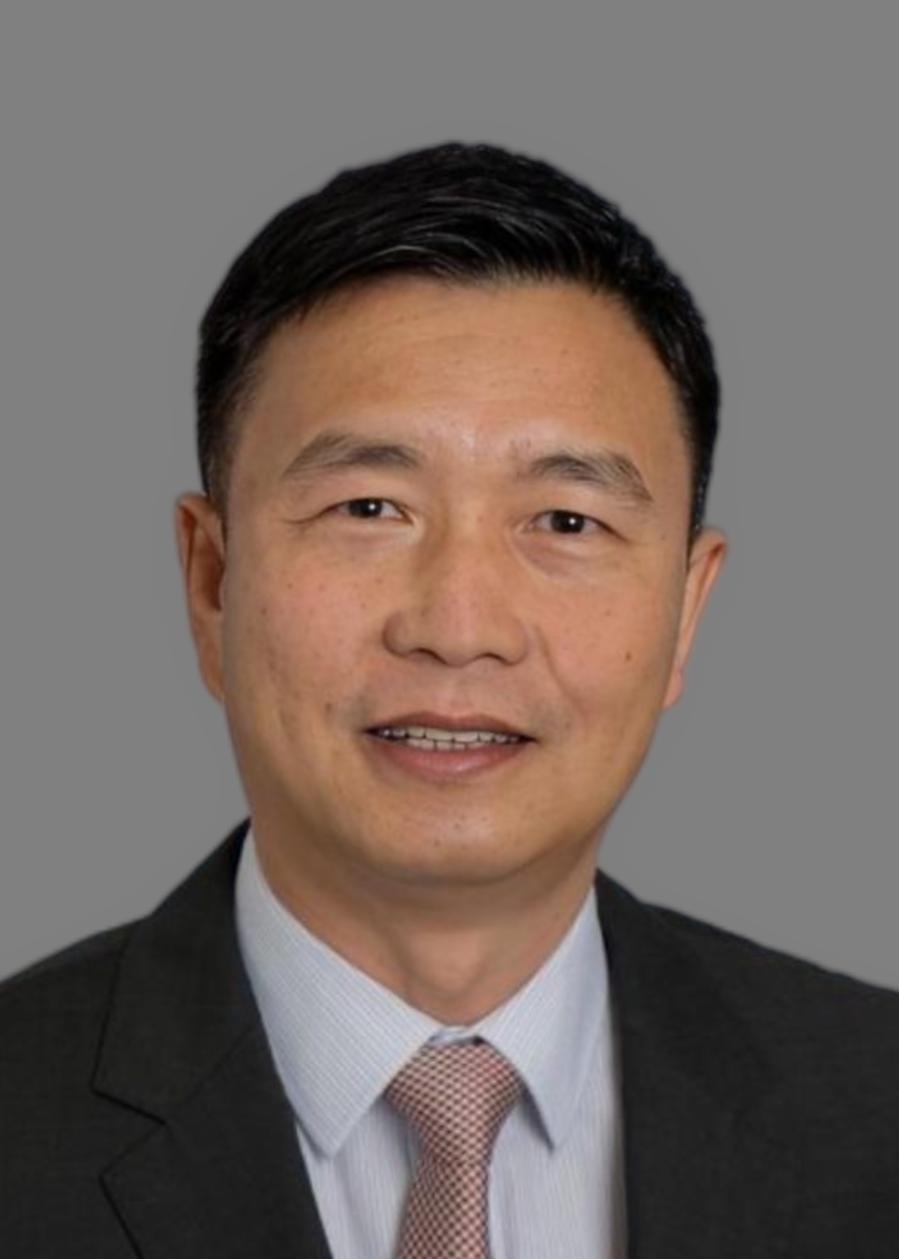}}]{Yong Xiang} (Senior Member, IEEE) is a Professor at the School of Information Technology, Deakin University, Australia. He received the B.E. and M.Eng. degrees in communication engineering from the University of Electronic Science and Technology of China, Chengdu, China, in 1983 and 1989, respectively. He received the PH.D. degree in computer science from the University of Melbourne, Parkville, Australia, in 2003. Dr. Xiang is also the Director of Trustworthy Intelligent Computing Lab (formerly named Deakin Blockchain Innovation Lab). He was the Associate Head of School (Research) (2013-2018), the Director of Artificial Intelligence and Data Analytics Research Cluster (2013-2019), the Director of Data to Intelligence Research Centre (2019-2020), and the Director of Deakin Univ.-Southwest Univ. Joint Research Centre on Big Data (2019-2021).
\end{IEEEbiography}

%








\end{document}